\documentclass[journal]{IEEEtran}
\usepackage{verbatim}
\usepackage{amsmath}
\usepackage{amsfonts}
\usepackage{amsthm}
\usepackage{tabls}
\usepackage{lipsum}
\usepackage{tabularx}
\usepackage{color}
\usepackage{graphicx}
\usepackage{amsmath}
\usepackage{pdflscape}
\usepackage{multicol}
\newcommand{\norm}[1]{\left\lVert#1\right\rVert}
\usepackage{amssymb}		
\usepackage{graphicx}		
\usepackage{hyperref}		

\usepackage{mathtools, cuted}
\usepackage{multirow}

\usepackage{booktabs}

\usepackage{caption}
\usepackage{subcaption}

\usepackage{algorithm,algcompatible,amsmath}
\algnewcommand\INPUT{\item[\textbf{Input:}]}%
\algnewcommand\OUTPUT{\item[\textbf{Output:}]}%

\floatstyle{plain}
\newfloat{twocolequfloat}{b}{zzz}
\floatname{twocolequfloat}{Equation}

\ifCLASSINFOpdf
\else
\fi

\begin{document}

\title{K-MACE and Kernel K-MACE  Clustering}
\author{ Soosan~Beheshti, Edward~Nidoy, and Faizan Rahman
\thanks{Authors are with the Department of Electrical, Computer, and Biomedical Engineering, Ryerson
University, Toronto, ON M5B 2K3, Canada. Email: \tt\small{ soosan@ee.ryerson.ca}}}

\maketitle

\begin{abstract}
Determining the correct number of clusters (CNC) is an important task in data clustering and has a critical effect on ﬁnalizing the partitioning results. K-means is one of the popular methods of clustering that requires CNC. Validity index methods use an additional optimization procedure to estimate the CNC for K-means. We propose an alternative validity index approach denoted by k-minimizing Average Central Error (KMACE). ACE is the average error between the true unavailable cluster center and the estimated cluster center for each sample data. Kernel K-MACE is kernel K-means that is equipped with an efficient CNC estimator. In addition, kernel K-MACE includes the ﬁrst automatically tuned procedure for choosing the Gaussian kernel parameters. Simulation results for both synthetic and real data show superiority of K-MACE and kernel K-MACE over the conventional clustering methods not only in CNC estimation but also in the partitioning procedure. 
\end{abstract}


\section{Introduction}
\label{intro}

\IEEEPARstart{C}{lustering} is one of the most used unsupervised learning tasks where unlabeled observed data samples are grouped based on their similarities and dissimilarities. It has vast applications in various areas such as image segmentation, market research, and sequence analysis. Two main challenges of clustering procedure are estimating the optimum number of clusters and partitioning the data \cite{mintro1},\cite{mintro2}. One of the well known partitional clustering methods is K-means \cite{kmeans}. Given the data, the method estimates cluster centers and partitions the data in a simple iterative procedure. Similar to other partitional approaches, K-means requires  the correct number of clusters (CNC) to finalize the clustering procedure. In general, most clustering algorithms require CNC estimate as their input \cite{ieee1},\cite{ieee2},\cite{ieee3},\cite{ref:ditc}. However, in many practical applications this value is not available and CNC has to be estimated during the clustering procedure by using the same data that requires clustering.
Overestimating the number of clusters typically results in redundant splitting  of a cluster,   and underestimating the number of clusters, on the other hand, results in combining clusters to form a loose clustering. 

Number of approaches for estimating the CNC in K-means are provided in \cite{kmcnc}. Pioneer methods of CNC estimation involve formulation of validity indexes. Validity indexes provide a quantitative measurement for comparison of clustering results with different number of clusters. To estimate CNC in K-means, validity indexes approaches use the results of $m$-clustering, e.i., clustering with $m$ clusters, for a range of $m$, $m \in [m_{min}, m_{max}]$. By solving an additional optimization problem the results of these $m$ clusterings are compared and the value of $m$ which optimizes the desired criterion is chosen as the estimate of CNC.  Examples of these validity indexes are gap index \cite{gap},  Silhouette index  \cite{sil}, Calinski-Harabasz index \cite{CB}, and Davies-Bouldin index\cite{DB}. More methods are proposed in \cite{vi1},\cite{vi2}. 
K-means is also used in divisive hierarchical clustering methods such as X-means \cite{xmeans}, G-means \cite{gmeans}, and dip-means  \cite{dip-means}. In these approaches stopping criterion for the hierarchical splitting procedure is a statistical test. X-means uses an information theoretical test criterion, G-means implements Anderson-Darling (AD) test and DIP-means, uses a criterion denoted by dip-dist. 

While clustering assignment on K-means is based on the distance of a sample to its cluster center, another family of clustering algorithms are density based where clusters are formed by grouping samples based on their proximity with respect to their neighboring samples. These methods provide the CNC estimate simultaneously. Density-Based Scan (DBSCAN) \cite{dbscan}, mean-shift \cite{mean_shift} and ordering points to identify cluster structure (OPTICS) \cite{optics} are some of the well known algorithms under this category (Note that OPTICS does not provide the clustering solution but rather provides a CNC estimate). There are other methods that concentrate on Fuzzy clustering or mixture modeling, \cite{fcm},\cite{BFC},\cite{MC},\cite{NML}. Concentration of this work is on providing a proper CNC estimator for K-means that overall keeps the computational complexity of clustering as low as K-means clustering itself. 

K-means uses the least square error to find the optimum cluster centers. On the other hand, the existing validity index methods, used with K-means, employ the available cluster compactness for CNC estimate to optimize a criterion that is not similar to the K-means partitioning criterion.  It seems rational to choose a consistent criterion for the validity index step as well. Motivated by this, MACE-means \cite{ref:9} proposed to minimize the Average Central Error (ACE).  It is worth mentioning that use of ACE in this problem setting shares the same fundamentals that use MSE in SVD order selection \cite{SVD}. However, MACE-means fails to estimate this error properly. It misses estimating the biased variance of this error and can lead to wrong CNC estimates. Here, we correct the definition of ACE to be the average least square error between the true center and the estimated center for each sample data in each $m$ clustering. Using the available cluster compactness, k-minimizing ACE (K-MACE) algorithm calculates estimate of both mean and variance of ACE to provide an accurate estimate of CNC \cite{ref:8}. While MACE-means was formulated for only uncorrelated spherical clusters, K-MACE is proposed for clustering both sphere and ellipsoid clusters. The proposed approach is very robust in estimating the covariance of each cluster.   On the other hand, kernel K-means clustering has been proposed for clustering overlapping and more arbitrary shaped clusters \cite{4}. In these scenarios kernel functions are used to transform data into a feature space. Similar to K-means the method requires estimation of CNC. Here, we extend K-MACE to kernel K-MACE to provide the CNC estimates along with the clustering procedure.  Note that in addition to the number of clusters $m$, existing Kernel based clustering methods require tuning the kernel function parameters. This  is currently done by trial and error and no method of validation and choosing the optimum parameters is available.  Consequently, the accuracy of the results depends on this ad hoc approach.  However, the Kernel parameter governs the separability of clusters in feature space and its optimum value corresponds to the true estimation of CNC. Here we propose another important feature  in the kernel K-MACE clustering algorithm that automatically tunes to the optimum Gaussian kernel parameters. While the computational complexity of K-MACE and kernel K-MACE are the same as that of K-means and kernel K-means, simulation results confirm advantages of the proposed approaches over the competing methods for both synthetic data and real data. It shows that K-MACE approaches outperform other methods even in the presences of major cluster overlaps.

The paper is organized as follows. Section II defines the considered clustering problem. 
 Average Central Error (ACE) and its importance in $m$-clustering evaluation are provided in Section III. Calculation of ACE, by only using the available data, for $m$-clustering and the K-MACE algorithm are in Section IV.  In Section V  kernel K-MACE is introduced. Section VI concentrates on simulation results for synthetic and real data sets and Section VII has the concluding remarks.

\section{Problem Statement}
Consider an observed data of length $N$, $x^N = [x(1) \; x(2)\; \dots \; x(N)]$ where $x(i) \in R_{d\times 1}$  represents collection of features. Each $x(i)$ is considered as a sample of a random vector with the following structure:
\begin{equation} \label{eq:1}
x(i) = \overline{c}_{x(i)} + \overline{w}_{x(i)}
\end{equation}
where scatter factor $\overline{w}_{x(i)}$ is a sample of $\overline{W}_{x(i)}$, a zero mean random vector that describes the variation of $x(i)$ around its center $\overline{c}_{x(i)}$.  For simplicity and without loss of generality, it can be assumed that the scatter factor has an additive white Gaussian distribution:
\begin{equation}\label{eq:dd_mod}
\overline{w}_{x(i)} \sim   \mathcal{N} (0, \overline{\Sigma}_{x(i)})
\end{equation}

where $\overline{\Sigma}_{x(i)}$ is a $d\times d$ covariance matrix. 

The notion of over bar is used for the true and unavailable elements of the cluster such as the true center and the true covariance of each cluster. The data set is generated by $\overline{m}$ mutually exclusive cluster model such that for the set of observed data we have $x^N \in X^N$:
\begin{equation}\label{memb}
X^N  =\overline{C}_1 \cup \overline{C}_2 \cup . . . \cup \overline{C}_{\overline{m}}
\end{equation}
Each cluster $\overline{C}(j)$ is described by its cluster center $\overline{c}_j$ paired with its covariance $\overline{\Sigma}(j)$. Therefore, for the whole data, we have:
\begin{eqnarray}\label{eq:added_1}
\rm{Cluster \ centers}&:& [\overline{c}(1) \; \overline{c}(2) \; \dots \; \overline{c}(\overline{m})] \\
\rm{Cluster \  Covariances}&:&[\overline{\Sigma}(1)\;  \overline{\Sigma}(2) \; \dots \; \overline{\Sigma}(\overline{m})] \nonumber
\end{eqnarray}

\begin{figure}
\centering
\includegraphics[scale=0.21]{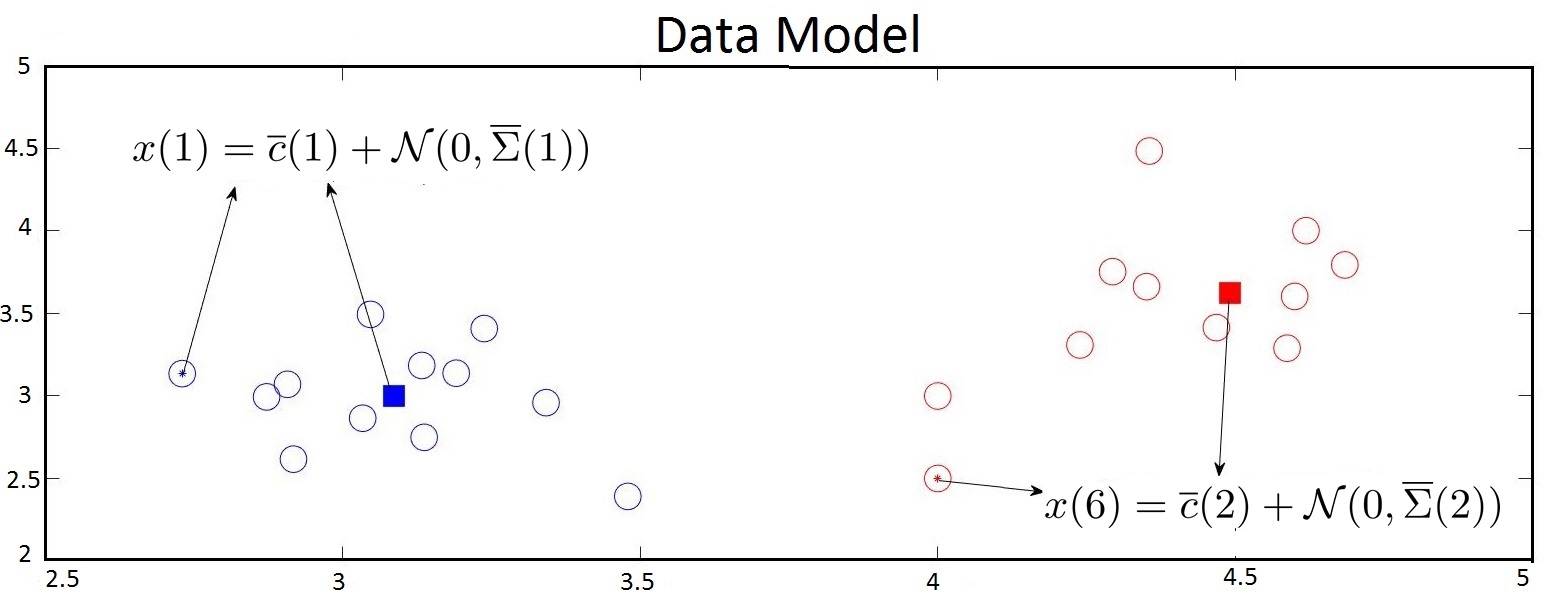}
\caption{An Example of Data Model in 2D. $(\overline{m}=2)$.}
\label{datamodel}
\end{figure}

Figure 1  shows an example of these clusters. Note that for each sample $x$ of observed data $X^N$ there is center of its true clusters. We denote its associated true center and true covariance with $\bar c_x$ and $\bar \Sigma_x$. For example as
Figure \ref{datamodel} shows the true center associated with $x(1)$ is the center of first cluster, $\overline{c}(1)$, while the true center associated with $x(6)$ is that of the second cluster, $\overline{c}(2)$. We rewrite this fact from the point of view of the samples $x(1)$ and $x(6)$ as follows:
\begin{equation}
 \overline{c}_{x(1)} = \overline{c}(1), \ \ \overline{\Sigma}_{x(1)} = \overline{\Sigma}(1)  \ \ \ \ \ \ \ \ \overline{c}_{x(6)} = \overline{c}(2), \ \   \overline{\Sigma}_{x(6)} = \overline{\Sigma}(2)
\end{equation}

The correct number of cluster (CNC) ($\overline{m}$ in (\ref{memb})) is not known. The data can be clustered with a range of possible cluster numbers $m$, ($ M_{min} \leq m \leq M_{max}$) where $M_{min}$ and $M_{max}$ are the considered upper and lower bounds for the true number of clusters. The goal is to provide a comparative measure to compare the results of $m$ clustering and choose the optimum $m$.

\section{Average Central Error (ACE) and $m$-clustering}
In this section we first formulate the $m$-clustering with proper notations. Next to compare the results of $m$ clustering, we suggest to use Average Central Error (ACE) as a criterion and choose the number of clusters, $m$, that minimizes this error.
\subsection{$m$-clustering Notations}
In $m$ clustering, the available data is used in finding $m$ clusters. Our method of clustering is K-means. The following are notations to distinguish for each data $x$ between its true unknown cluster and its $m$-clustering membership.

Clustering the available data into $m$ mutually exclusive clusters results in $m$ centers. Each of these $m$ cluster is denoted by $C_{mj}$ where the subscript $j$ pertains to the $j^{th}$ cluster. Similar to (\ref{eq:added_1}) the centers of this clustering are

\begin{eqnarray}\label{estt}
\rm{Estimated \ centers}&:& [{\hat c}_{m1} \; \cdots \; {\hat c}_{mj} \; \dots \; {\hat c} _{mm}] 
\end{eqnarray}

Note that the notion of over bar is now replaced by hat to represent the estimation and the new subscript has two elements to show the number of clusters $m$ used in clustering as well as indexing the cluster with $j$.  

 In this case the available data $x$ in (1) is now denoted by $x_{mj}$ which means that the data has been clustered by $m$ clusters and is now member of the $j$th cluster. Each  element of cluster with center $c_{mj}$, that is now denoted by $x_{mj}$,  already has been generated by a true unknown center.  In each $m$ clustering, each $x(i)$ becomes a $x_{mj}$, its true center and covariance is now denoted by $\overline{c}_{x_{mj}}$ and we can represent (1) in the following form: 

\begin{equation} \label{eq:1}
x_{mj} = \overline{c}_{x_{mj}} + \overline{w}_{x_{mj}}
\end{equation}

Therefore each element of $j$th cluster in $m$ clustering has been generated by its true center and is pointing to the center of this cluster:
 
\begin{equation}\label{erl}
x_{mj} \in C_{mj}:\;\; \overline{c}_{x_{mj}} \rightarrow x_{mj} \rightarrow \hat{c}_{mj} 
\end{equation}
where the estimated center of the $j^{th}$ cluster is denoted by $\hat{c}_{mj}$, and is calculated by averaging the $j^{th}$ cluster members: 
\begin{equation}\label{cmj}
\hat{c}_{mj} = \frac{1}{n_{mj}} \sum_{i=1}^{n_{mj}} x_{mj}(i)
\end{equation}
and $n_{mj}$ is the number of elements in $C_{mj}$.

Figure \ref{fig:clust_example} illustrates an example of  $m$ clustering. As the figure shows the true number of clusters is two, $\overline{m}=2$, with two true cluster centers $\bar c_1$ and $\bar c_2$. The figure shows the clustering results for $m=1,2,3$. The figure also shows one element of the first cluster as $x$. In the case of 2-clustering ($m=2$), 
 Figure \ref{fig:clust_example}(a) shows the two estimated centers $\hat{c}_{21}$ and $\hat{c}_{22}$. Note that in this scenario $x$ is also a member of $x_{21}$ as it belongs to the first estimated cluster. On the other hand, in the case of $m=3$, in  Figure \ref{fig:clust_example}(c), $x$ is a member of the second estimated cluster $C_{32}$ and is a member of set $x_{32}$. 

\subsection{ACE in m-clustering}
ACE, denoted by $Z_{m}$ is a measure of error between estimated cluster centers in (\ref{estt}) and true cluster centers defined in (\ref{eq:added_1}). In each $m$-clustring, for the $j^{th}$ cluster $C_{mj}$, this error is
\begin{equation}\label{zsmj}
\rm{Center \ Error \ for \ cluster \ } C_{mj}: z_{mj} = \sum_{i=1}^{n_{mj}} \norm{ \bar{c}_{x_{mj}(i)} -\hat{c}_{mj} }_2^2 
\end{equation}
Note that the summation is over elements of each cluster.

The ACE for $m$ clustering is defined as summation of all the cluster center errors in (\ref{zsmj}) for all $m$ clusters divided by the total number of element $N$:

\begin{equation} \label{eq:6}
z_{m}  = \frac{1}{N} \sum_{j=1}^{m}z_{{mj}}
\end{equation}
\begin{figure*}
\centering
\begin{subfigure}{0.4\textwidth}
  \raggedright
  \includegraphics[width=1.1\linewidth]{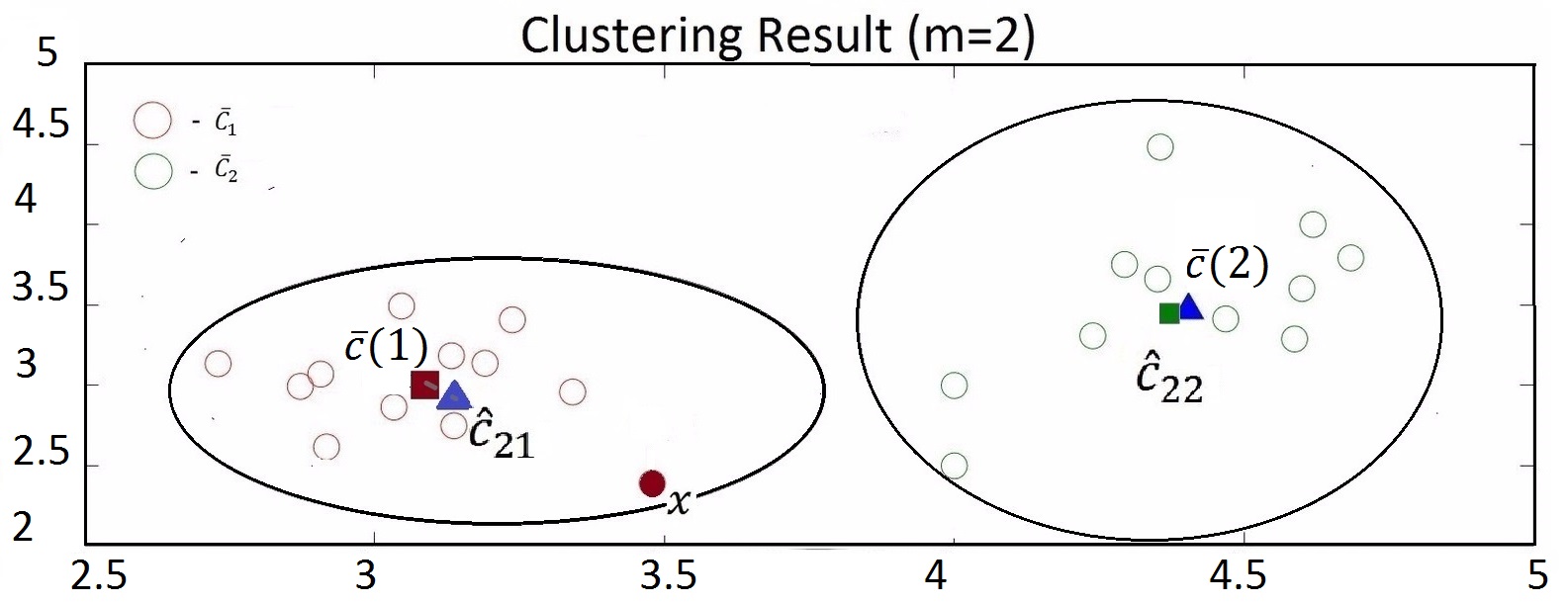}
  \caption{$m$ Clustering $(m=2)$ and $x \in C_{m1}$}
   \label{overestimated_m}
\end{subfigure}\hspace{15mm}

\begin{subfigure}{0.4\textwidth}
  \raggedleft
  \includegraphics[width=1 \linewidth]{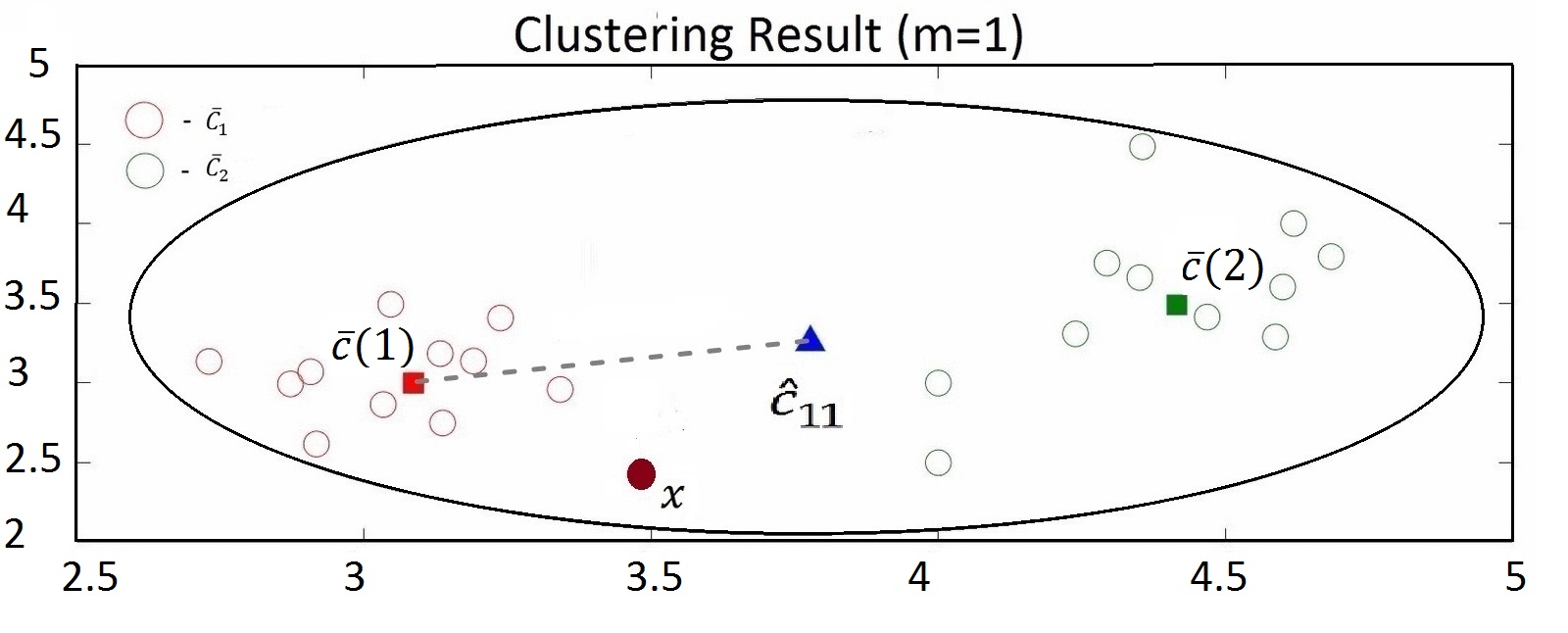}
  \caption{$m$ Clustering $(m=1)$ and $x \in C_{m1}$}
  \label{underestimated_m}
\end{subfigure}
\begin{subfigure}{0.4\textwidth}
  \raggedright
  \includegraphics[width=1\linewidth]{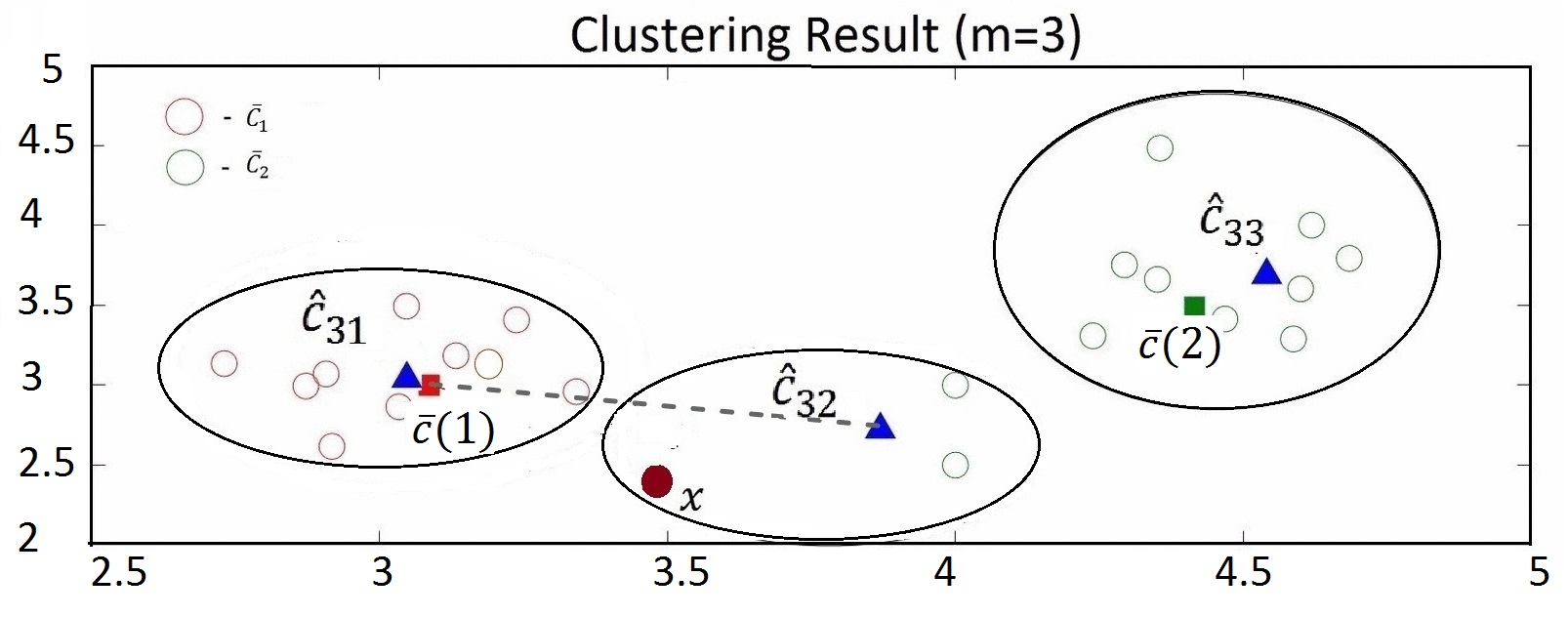}
  \caption{$m$ Clustering $(m=3)$ and $x \in C_{m2}$}
   \label{overestimated_m}
\end{subfigure}
\caption{Example of Central Error for observed sample $x$.}
\label{fig:clust_example}
\end{figure*}

\begin{figure*}
	\centering
  	\includegraphics[width=0.8\linewidth]{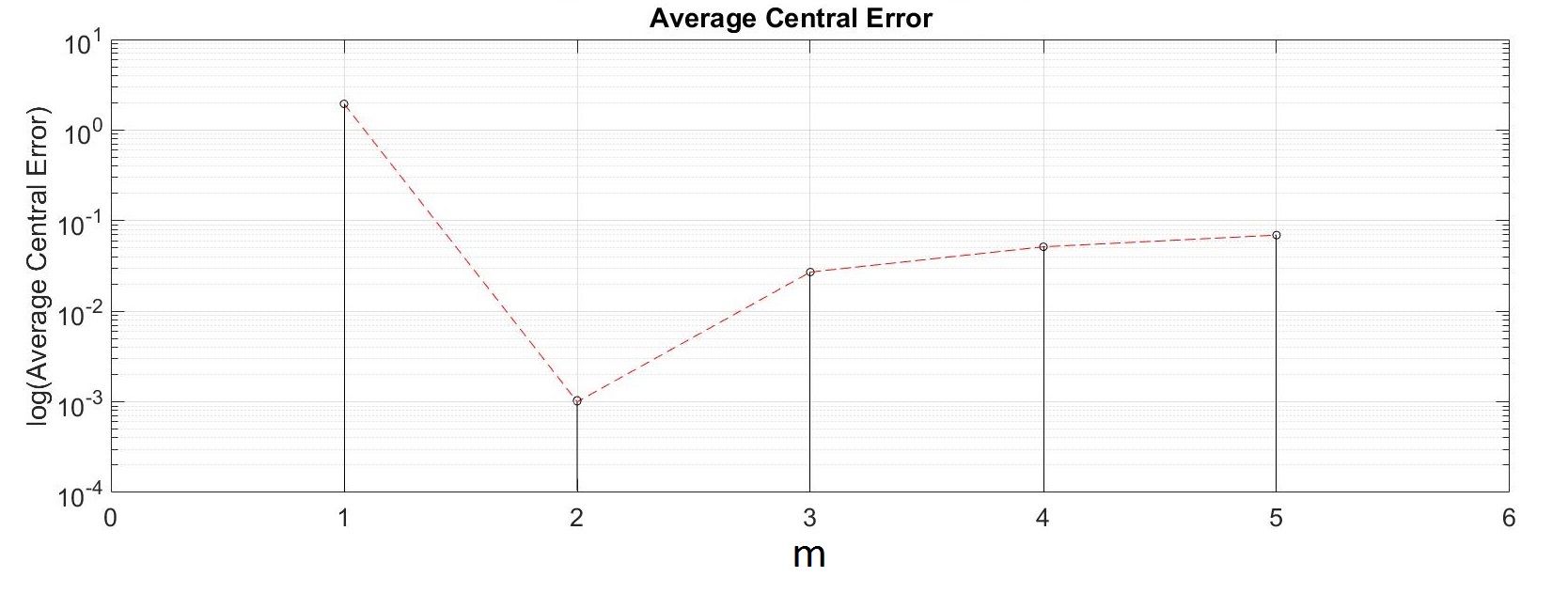}
  	\caption{Average Central Error as a Function of $m$}
   	\label{minace}
\end{figure*}

Figure \ref{fig:clust_example} also illustrates an example of behavior of the central error in a scenario that the true number of clusters is two, $\overline{m}=2$. As $x$  is generated by the first cluster, its associated correct center is $\bar {c}(1)$.  Therefore, the central error for the shown $x$ in the case of two clustering ($m=2$) is $\norm{\overline{c}(1) - \hat{c}_{21}}_2^2$. In the case of $m=1$, the central error associated for $x$ is $ \norm{\overline{c}(1) - \hat{c}_{11}}_2^2$ which is a much larger value. On the other hand, in the case of $m=3$  the associate central error for $x$ is  $\norm{\overline{c}(1) - \hat{c}_{32}}_2^2$ which is again larger than the first case ($m=2$).  The associate central error for the other data has a similar behavior.  Therefore, as we perform $m$ clustering where $m$ becomes further(smaller/larger) from the true number of cluster $\overline{m}$, we expect this distance to increase for each $x$. Consequently, adding all these central errors for all the data,  in (\ref{eq:6}), expected to have the minimum value at $m$ clustering when $m=\overline{m}$. Figure \ref{minace} shows behavior of average central error with respect to $m$ for the data set depicted on Figure \ref{fig:clust_example}.

In the following sections we study ACE closely and provide a method of estimating the unavailable ACE by $only$ using the available  data and without the knowledge of the true cluster centers.

\subsection{Mean and Variance of ACE}
The average central error, $z_{{mj}}$ in (\ref{zsmj}), is a sample of random variable $Z_{{mj}}$. Here, we provide the expected value and variance of ACE in (\ref{eq:6}).

\textit{Lemma 1:} The central error for the $j^{th}$ cluster in $m$ clustering, $z_{{mj}}$ in (\ref{zsmj}), is a sample of random variable $Z_{{mj}}$ that is a function of random vector $\overline{W}_{x_{mj}(i)}$, defined in (\ref{eq:dd_mod}):
\begin{equation}
\label{eq:10}
\begin{split}
Z_{{mj}}  = \norm{\Delta_{{mj}}}_2^2+ \frac{1}{n_{mj}}\sum_{i=1}^{n_{mj}}\overline{W}_{x_{mj}(i)}^T\overline{W}_{x_{mj}(i)} + \\ \frac{1}{n_{mj}}\sum\limits_{i \neq k}^{n_{mj}}\overline{W}_{x_{mj}(i)}^T\overline{W}_{x_{mj}(k)}^T 
\end{split}
\end{equation}
with the following expected value and variance of $Z_{{mj}}$  

\begin{equation}\label{eq:12}
E[Z_{{mj}}]  = \norm{\Delta_{{mj}}}_2^2  + \frac{1}{n_{mj}}\sum_{i=1}^{n_{mj}}tr(\overline{\wedge}_{x_{mj}(i)})
\end{equation}
\begin{equation}\label{eq:13}
\begin{split}
Var[Z_{{mj}}]  = \frac{2}{{n_{mj}}^2}\sum_{i=1}^{n_{mj}}tr((\overline{\wedge}_{x_{mj}(i)})^2) + \\
 \frac{2}{{n_{mj}}^2}\sum\limits_{i \neq k}^{n_{mj}}tr(\overline{\wedge}_{x_{mj}(i)}\overline{\wedge}_{x_{mj}(k)})
 \end{split}
\end{equation}
where $tr()$ is the trace function of a diagonal matrix, $\overline{\wedge}_{x_{mj}(i)}$ whose values are the eigenvalues of $\overline{W}_{x_{mj}(i)}$'s 
covariance  and $\norm{\Delta_{{mj}}}_2^2$ is
\begin{eqnarray}
 \norm{\Delta_{{mj}}}_2^2 = \norm{\overline{c}_{x_{mj}}A_{mj}}_2^2,
 \end{eqnarray}
Matrix $A_{mj}$ is a $n_{mj}$ by $n_{mj}$ matrix   
\begin{equation}\label{eq:11}
A_{mj} = 
\begin{bmatrix}
\frac{n_{mj} - 1}{n_{mj}}&  \frac{-1}{n_{mj}} & \dots{} & \frac{-1}{n_{mj}} \\
\frac{-1}{n_{mj}} &  \frac{n_{mj} - 1}{n_{mj}} &\frac{-1}{n_{mj}} &\vdots\\
\vdots{}&\dots{}&\ddots{} & \frac{-1}{n_{mj}}\\
\frac{-1}{n_{mj}} &\dots{} & \frac{-1}{n_{mj}} & \frac{n_{mj} - 1}{n_{mj}}
\end{bmatrix}
\end{equation}
and $\bar c_{x_{mj}}$ is a vector of all the associated $\bar c_{x_{mj(i)}}$s for the elements of $C_{mj}$ .

\textit{Proof} in Appendix \ref{app_lemma1}.

Consequently, the expected value and variance of the overall ACE, $Z_{m}$, in (\ref{eq:6}), are:

\begin{equation}\label{eq:14}
E[Z_{m}]  = \frac{1}{N} \sum_{j=1}^m E[Z_{{mj}}]
\end{equation}
\begin{equation}\label{eq:15}
Var[Z_{m}]  = \frac{1}{N^2} \sum_{j=1}^m Var[Z_{{mj}}]
\end{equation}

Note that the variance of $Z_{m}$ is also summation of the individual variance of $Z_{{mj}}$'s in (\ref{eq:15}). Due to the independence of $W_{x(i)}$ and $W_{x(k)}$ $(i \neq k)$, $E[Z_{{mj}}Z_{{mi}}]= E[Z_{{mj}}]E[Z_{{mi}}]$ and therefore, the covariances between $Z_{{mj}}$ and $Z_{{mi}}$ is equal to zero   for $j \neq i$.

From (\ref{eq:12}), the expected value of $Z_{{mj}}$ has two terms. The first term, $\norm{\Delta_{{mj}}}_2^2 $, is a function of the unknown true cluster center. This term is a decreasing function of $m$. The second term is a function of the cluster covariance and is monotonically proportional to $m$. As a result, there is a point in which $E[Z_{m}]$ reaches its minimum value at some $m$. This is a manifestation of a form of  bias-variance trade-off \cite{ref:8},\cite{ref:8.1}.

\subsection{ACE Mean Estimation}
Given the available data, no direct information on ACE, $z_{m}$, is available. Here we propose a method to use the available data to estimate ACE. 

\subsubsection{Cluster Compactness}
The Cluster Compactness denoted by $y_{m}$ is the available error between samples of a cluster and its  estimated cluster center. The Cluster Compactness of the $j^{th}$ cluster in $m$-clustering is 
\begin{equation}\label{eq:9}
y_{{mj}}=\sum_{i=1}^{n_{mj}} \norm{{{x_{mj}(i)} - \hat{c}_{mj}}}_2^2
\end{equation}

\begin{equation}\label{eq:8}
m\rm{-clustering \ Compactness: \ }\mathit{y_{m}  = \frac{1}{N} \sum_{j=1}^{m}y_{{mj}}}
\end{equation}
To estimate ACE we use the available cluster compactness. The following Lemma connects the structure of mean and variance of cluster compactness with that of ACE.

\textit{Lemma 2}: The observed cluster compactness $y_{{mj}}$, in (\ref{eq:9}), is a sample of random variable $Y_{{mj}}$ with the following expected value and variance:
\begin{equation} \label{eq:23}
E[Y_{{mj}}]  =  \norm{\Delta_{{mj}}}_2^2  + \frac{(n_{mj} - 1)}{n_{mj}} \sum_{i=1}^{n_{mj}}tr(\overline{\wedge}_{x_{mj}(i)})
\end{equation}

\begin{equation}\label{eq:24}
\begin{split}
Var[Y_{{mj}}]  = \frac{2(n_{mj}-1)^2}{n_{mj}^2}\sum_{i=1}^{n_{mj}} tr((\overline{\wedge}_{x_{mj}(i)})^2) + \ \ \ \ \ \ \ \ \ \ \ \ \ \ \ \ \ \ \ \ \ \ \ \ \ \ \ \ \ \\
\frac{1}{n_{mj}^2}\sum_{i \neq k}^{n_{mj}}tr(\overline{\wedge}_{x_{mj}(i)} \overline{\wedge}_{x_{mj}(k)}) + \frac{4}{d\times n_{mj}}   \norm{\Delta_{{mj}}}_2^2   \sum_{i=1}^{n_{mj}} tr(\overline{\wedge}_{x_{mj}(i)}) \ \ \ \ \ \ \ \ \ \ 
\end{split}
\end{equation}

\textit{Proof} in Appendix \ref{app_lemma_2}.

Consequently, the expected value and variance of overall $m$-clustering compactness, $Y_{m}$ are: 
\begin{equation}\label{eq:25}
E[Y_{m}]  = \frac{1}{N} \sum_{j=1}^m E[Y_{{mj}}]
\end{equation}
\begin{equation}\label{eq:26}
Var[Y_{m}]  = \frac{1}{N^2} \sum_{j=1}^m Var[Y_{{mj}}]
\end{equation}
Note that the variance of $Y_{m}$ is a summation of individual variances of $Y_{{mj}}$ due to the independence of $W_{x(i)}$ and $W_{x(k)}$ $(i \neq k)$.
\subsubsection{Probabilistic bounds of ${\norm{\Delta_{{mj}}}_2^2 }$ using cluster compactness}

Inspecting \textit{Lemma 1} and \textit{Lemma 2}, it can be seen that the mean of both $Z_{{mj}}$ and $Y_{{mj}}$ are functions of the unknown term $\norm{\Delta_{{mj}}}_2^2$. The following theorem provides probabilistic bounds on this value by using the available $y_{{mj}}$. \\

\textit{Theorem 1}: With validation probability  $P_v=1-\frac{1}{\alpha_{mj}^2}$, $\norm{\Delta_{{mj}}}_2^2$ is bounded by 

\begin{equation}
\label{bound_delta}
\underline{\norm{\Delta_{{mj}}}_2^2}  \leq \norm{\Delta_{{mj}}}_2^2 \leq \overline{\norm{\Delta_{{mj}}}_2^2} 
\end{equation}
where 

\begin{equation}\label{delta_upp}
\begin{split}
\overline{\norm{\Delta_{{mj}}}_2^2} = \bar{y}_{{mj}} + k_{{mj}} \\
 \underline{\norm{\Delta_{{mj}}}_2^2} = \bar{y}_{{mj}} - k_{{mj}}
\end{split}
\end{equation}

and

\begin{equation}
\bar{y}_{{mj}} = (g_{mj}- y_{{mj}}) -  \alpha_{mj}^2 \frac{2}{d\times n_{mj}}\sum_{i=1}^{n_{mj}} tr(\overline{\wedge}_{x_{mj}(i)} )
\end{equation}
where $g_{mj}= \frac{(n_{mj} - 1)}{n_{mj}}\sum_{i=1}^{n_{mj}}tr(\overline{\wedge}_{x_{mj}(i)}) $, and
\begin{align}
\label{radical}
k_{{mj}} =& 2\alpha_{mj}  [ \nonumber\\
&v_{{mj}}+ \left(  \frac{(4\alpha_{mj}^2+2d^2)}{d^2n_{mj}^2}\sum_{i \neq k}^{n_{mj}}tr(\overline{\wedge}_{x_{mj}(i)}  \overline{\wedge}_{x_{mj}(k)})  \right) \nonumber \\ &+
 \left(  \frac{ (4\alpha_{mj}^2+2d^2 n_{mj}^2) } {d^2n_{mj}^2} \sum_{i=1}^{n_{mj}} tr((\overline{\wedge}_{x_{mj}(i)})^2) \right)   ] ^{1/2} 
\end{align}
and
\begin{equation}
v_{{mj}} = \frac{-4( g_{mj}- y_{{mj}} )} {d\times n_{mj}}\sum_{i=1}^{n_{mj}}tr(\overline{\wedge}_{x_{mj}(i)})
\end{equation}

\textit{Proof}: Using Chebyshev's inequality  
 
\begin{equation}\label{eq:30}
P(|E[Y_{{mj}}] - y_{{mj}}| \leq \alpha_{mj} \sqrt{var[Y_{{mj}}]}) \geq P_v,
\end{equation} 
with $P_v=1-\frac{1}{\alpha_{mj}^2}$, the observed $y_{{mj}}$ is bounded by 
\begin{equation}\label{ysmj_bounds}
E[Y_{{mj}}] - \alpha_{mj} \sqrt{var[Y_{{mj}}]} \leq y_{{mj}} \leq E[Y_{{mj}}] + \alpha_{mj} \sqrt{var[Y_{{mj}}]}
\end{equation}
Using (\ref{eq:23}) and (\ref{eq:24}) in (\ref{ysmj_bounds}) and solving for the boundaries of the two pairs of inequality in (\ref{ysmj_bounds}) results in bounds $\overline{\norm{\Delta_{{mj}}}_2^2}$ and $\underline{\norm{\Delta_{{mj}}}_2^2}$. Details are in Appendix C.

\section{K-MACE}
Using the data to calculate bounds on $\norm{\Delta_{{mj}}}_2^2$, we have bounds on the  mean of ACE in (\ref{eq:12}). Here we provide bounds for ACE for $m$-clustering and suggest to use the $m$-clustering that minimizes this criterion.   
\subsection{Bounds on ACE}
Employing Chebyshev's inequality, with probability $P_c=1-\frac{1}{\beta^2}$, the ACE, $z_{m}$, lies around the expected value of random variable $Z_{m}$ such that:
\begin{equation}\label{eq:16}
P(|E[Z_{m}] - z_{m}| \leq \beta \sqrt{var[Z_{m}]}) \geq P_c
\end{equation}

Therefore, with confidence probability $P_c$, $z_{m}$ is bounded as follows:
\begin{equation}\label{eq:19}
\underline{z_{{m}}} \leq z_{{m}} \leq \overline{z_{{m}}}
\end{equation} 
where $\underline{z_{{m}}}$ and $\overline{z_{{m}}}$ are the resulted lower bound and upper bound of $z_{m}$ respectively
\begin{equation}\label{eq:18}
\overline{z_{{m}}} =  E[Z_{{m}}] + \beta \sqrt{var[Z_{{m}}]}
\end{equation} 
\begin{equation}\label{eq:17}
\underline{z_{{m}}} =  E[Z_{{m}}] - \beta \sqrt{var[Z_{{m}}]}
\end{equation} 
In this calculation $E[Z_{{m}}]$ in (\ref{eq:12}) is replaced by its estimate using bounds on $\norm{\Delta_{{mj}}}_2^2$ in (\ref{bound_delta}).
\subsection{Validation and Confidence Probabilities}\label{vall}
Values of $P_v$ and $P_c$ validation and confidence probabilities are close to one as long as $\alpha_{mj}$s and $\beta$ are chosen large enough and definitely larger than one. Simultaneously, the upper bounds and lower bounds on $\norm{\Delta_{{mj}}}_2^2$ in (\ref{bound_delta}) and $z_{m}$ in (\ref{eq:19}) are close to each other as long as $\frac{\beta}{N}$ and 
$\frac{\alpha_{mj}}{\sqrt{d \times n_{mj}}}$ are small.  Note that both cluster compactness $Y_{{mj}}$, in (\ref{eq:9}) and ACE, $z_{{mj}}$ in (\ref{zsmj}), are summations of squared Euclidean distances. The number of elements in these summations is the number of elements of each cluster and when this number is large enough, by using the Central Limit Theorem (CLT), the inequalities in (\ref{eq:16}) and (\ref{eq:30}) can be turned into equalities in the following forms with $P_v=Q(\alpha_{mj})$, $P_c=Q(\beta)$:
\begin{eqnarray}
P(|E[Y_{{m}}] - y_{{m}}| \leq \alpha_{mj} \sqrt{var[Y_{{m}}]}) = Q(\alpha_{mj})\\
P(|E[Z_{{mj}}] - z_{{mj}}| \leq \beta \sqrt{var[Z_{{mj}}]}) = Q(\beta) 
\end{eqnarray}

For example when the scatter factor $W$ has a Gaussian structure, these two errors are Chi-squared and in practical applications sum of chi-square when the number of elements are ten or more is well approximated with a Gaussian distribution.   

\subsection{K-MACE Using the Available Data}
\label{finalstep}
When clustering the data with $k$ number of clusters,  associated with each sample $x(i)$ that now belongs to cluster $C_{kj}$, the covariance estimate of scatter factor in (\ref{eq:dd_mod}) is the covariance of cluster $C_{kj}$. We denote this covariance estimate for each $x_{kj}$ as $\hat \Sigma_{x_{kj}}$. This covariance estimate is then used in calculating bound of $z_{m}$ in (\ref{eq:18}) and (\ref{eq:17}) , with $m$-clustering for all values of $m$. The estimated $z_{m}$ with this covariance estimate is denoted by $z_{{m,k}}$ and calculated upper and lower bounds for the ACE in this case are denoted by $\overline{z_{{m,k}}}$ and $\underline{z_{{m,k}}}$. Minimizing this upperbound, we have
\begin{equation}
\label{mk}
\hat{m}_k = \arg \min_m (\overline{z_{{m,k}}})
\end{equation}
If $k$ is the CNC, it is expected for $\hat{m}_k$ to be equal to $k$ itself. If this condition does not hold, it means that $k$ is further from the CNC. Consequently, the optimum CNC estimate $\hat{m}$ is the $k$ value for which this normalized
 discrepancy is minimum 
\begin{equation}
\label{hhhhh}
k^* = \arg \min_k ( \frac{\overline{z_{{\hat{m}_{k},k}}} - \overline{z_{{k,k}}}} {\overline{z_{{\hat{m}_{k},k}}}})
\end{equation}

and 

\begin{equation}
\label{finalm}
\hat{m} = \hat{m}_{k^*}
\end{equation}

Detailed pseudo-code for implementation of K-MACE algorithm is described in Algorithm \ref{kmace_alg}. Iin this procedure, probabilistic bounds on $||\Delta_{mj}||_2^2$ are calculated first in order to provide ACE bounds. Note that in line 10 of the algorithm if this term have a negative estimated value, it indicates that we are very far from the true value of $k$ and therefore the case should be discarded. Consequently, for these values of $k$ we set the upperbound to a large values $\Delta_{max}$ to get excluded form the comparison. As the algorithm shows, the complexity of K-MACE method in finding the CNC is linear, $O(dN)$. 
 
\begin{algorithm*}
\centering
   \caption{K-MACE Algorithm}
  \begin{algorithmic}[1]
    \INPUT Data set $\textbf{x} = [x_1,x_2,...,x_N]$, range of m, $[m_{min},m_{max}]$
    \OUTPUT Estimate of the number of cluster $\hat{m}$ and optimum clustering solution.
    \FOR{($m=m_{min}; m \leq m_{max};m_{++}$)}
        \STATE $[C_{m1},C_{m2},..,C_{mm}] =kmeans(x,m)$
         \FOR{each cluster $C_{mj} \  j=1,..,m$}     	 
            \STATE Solve cluster compactness $y_{mj}$ of cluster $C_{mj}$ using (\ref{eq:9})
       \ENDFOR   
       \STATE Solve for total cluster compactness $y_{m}$ using (\ref{eq:8}) 
   \ENDFOR 

    \FOR{($k=m_{min}; k \leq m_{max};k_{++}$)}
    	\STATE Calculate for cluster covariances denoted by for $\hat{\wedge}_{x_i}$ from clustering solution $[C_{k1},C_{k2},..,C_{kk}]$ $\hat{\wedge}_{x_i}$ =  eigenvalues of the covariance matrix of the cluster which $x_i$ belongs to.
        \STATE  Use $\hat{\wedge}_{x_i} \rightarrow \overline{\wedge}_{x_i}$ to calculate for the upperbound of $\norm{\Delta_{{mj}}}_2^2$, $\overline{\norm{\Delta_{{mj}}}_2^2}$ using equations (\ref{delta_upp})-(\ref{radical}). In case the value is negative, set $\overline{\norm{\Delta_{{mj}}}_2^2}$ to $\Delta_{max}$.
    \STATE  From the results of steps 9 and 10, calculate for $E[Z_{{mj}}]$ using (\ref{eq:12}) and $Var[Z_{{mj}}]$ using (\ref{eq:13})
    \STATE  Calculate for $E[Z_{m}]$ and $Var[Z_{m}]$ in (\ref{eq:14}) and (\ref{eq:15}).
	\STATE The upperbound of $z_{m}$ under the assumption that there are $k$ cluster covariances, $\overline{z_{{m,k}}}$, can then be found using (\ref{eq:18}).
     \STATE From $\overline{z_{{m,k}}}$, we can obtain $\hat{m}_k$ using (\ref{mk}) 
     \STATE Save values $\overline{z_{{\hat{m}_k,k}}}$ and  $\overline{z_{{k,k}}}$.
   \ENDFOR
   
   \STATE From values saved on Step 15, the estimated CNC $\hat{m}$ can be found using (\ref{hhhhh}) and then (\ref{finalm}) and the optimum clustering solution is given by $[C_{\hat{m}1},C_{\hat{m}2},..,C_{\hat{m}\hat{m}}]$(One of the clustering solution from Step 2).
   
  \end{algorithmic}
  \label{kmace_alg}
\end{algorithm*}

\section{Kernel K-MACE}
In this section we utilize the K-MACE algorithms for the cases of kernel K-means to find the optimum number of clusters in this scenario. In addition, by using the kernel K-MACE we are also able to provide the optimum tuning parameter of the Gaussian kernels. 
\subsection{Initial Cluster Assignment in Kernel K-MACE}
The kernel K-means algorithm initially randomly assigns data samples to clusters. This is similar to K-mean, however, kernel k-means algorithm displays increased sensitivity to random assignment of data samples to clusters which decreases the consistency of the clustering results. Here we expand on the proposed K-means technique in \cite{mod_kmeans_p} for the kernel k-means scenario to eliminates the sensitivity problem. The approach uses the distance between data samples in feature space to iteratively assign data samples to clusters hence reducing the number of iterations required for the clustering algorithm to reach convergence.

\subsection{K-MACE in feature space}
Calculation of the optimum cluster and CNC is analogous to that of the K-MACE itself. In the feature space we have the following replacement of data $x$ with its feature counterparts $\phi$: 
\begin{equation} \label{link}
x_{mj} \to \phi_{mj} 
\end{equation}
Consequently $\overline{c}(\phi_{mj})$ generates the data and $ \hat{c}_{\phi_{mj}}$ is the respective cluster center estimate in $m$-clustering.  ACE and cluster compactness defined in (\ref{eq:6}) and (\ref{eq:8}) are now respectively in the following forms:
\begin{equation} \label{eq:6000}
\begin{gathered}
Z_{s_m}  = \frac{1}{N} \sum_{j=1}^{m}Z_{s_{mj}},\;\; Z_{s_{mj}} = {\lVert \overline{c}_{\phi_{mj}} - \hat{c}_{\phi_{mj}} \rVert }_F^2
\end{gathered}
\end{equation}
\begin{equation}\label{eq:18000}
\begin{gathered}
Y_{s_m}  = \frac{1}{N} \sum_{j=1}^{m}Y_{s_{mj}},\;\; Y_{s_{mj}} = {\lVert \hat{c}_{\phi_{mj}} - {\phi}_{mj} \rVert}_F^2
\end{gathered}
\end{equation}

\subsection{Optimum Gaussian Kernel Parameter}

The Gaussian kernel function parameter $\sigma$ dictates the structure of data in feature space and therefore the clustering of a dataset.  Very small or very large values of $\sigma$ result in data losing its structure which makes correctly clustering a dataset very difficult. Clustering results are calculated for a range of $\sigma_i$ : $[\sigma_{k_{min}},...,\sigma_{k_{max}}]$ values. 
To obtain the best $\sigma_i$, we propose to obtain the gradient of the minimum $\overline{Z_{s_{m}}}$ at each $\sigma_i$ and for each $\sigma_i$ after the value is  peaked, we add the absolute value of the gradient between the previous and the next value of $\sigma_i$. The value of $\sigma_i$ which corresponds to the maximum value of this sum is chosen and the corresponding $\hat{m}$ and clustering result are chosen as the correct result \cite{kkmace_p}:

{\scriptsize\begin{equation}\label{sigma_grad}
\begin{aligned}
\hat \sigma= \;\;\;\;\;\;\;\;\;\;\;\;\;\;\;\;\;\;\;\;\;\;\;\;\;\;\;\;\;\;\;\;\;\;\;\;\;\;\;\;\;\;\;\;\;\;\;\;\;\;\;\;\;\;\;\;\;\;\;\;\;\;\;\;\;\;\ \;\;\;\\
\max\limits_{\sigma_{k}}(|\frac{\partial (\min(\overline{Z_{s_m}}(m-1 \to m,\sigma_k)))}{\partial \sigma_k} + \frac{\partial (\min(\overline{Z_{s_m}}(m \to m+1,\sigma_k)))}{\partial \sigma_k}|) \\   
for   \quad   \sigma_k > \max\limits_{\sigma_{k}}(\min_m (\overline{Z_{s_m}}(m,\sigma_k)))
\end{aligned}
\end{equation}}
The minimum $\overline{Z_{s_{m}}}$ rises for increasing $\sigma_i$ as the clusters do not have any structure for very small values of $\sigma_i$ i.e. for $\sigma_i$ $<$ 1. As the value of $\sigma_i$ increases to the optimum, a rapid decrease in minimum $\overline{Z_{s_{m}}}$ occurs corresponding to the estimated value of $\sigma_i$ and the correct clustering result. 

Algorithm 2 is a pseudo code of Kernel K-MACE and line 2 can be used for Gaussian Kernel K-MACE to tune the Gaussian parameter. Computational complexity of kernel K-MACE remains the same as kernel K-means which is $O(N^2)$.



\begin{algorithm}[ht!]

\small

\centering

\caption{Kernel k-MACE Algorithm}

\begin{algorithmic}[1]

\INPUT Data set $\textbf{x} = [x_1,x_2,...,x_N]$, range of m
$[m_{min},m_{max}]$ and range of values for $\sigma_k$ =
$[\sigma_{k_{min}},...,\sigma_{k_{max}}]$ with a fixed interval

\OUTPUT Estimated number of clusters $\hat{m}$, the clustering solution
$[\hat{C}_1,\hat{C}_2,..,\hat{C}_{\hat{m}}]$ and the optimum value of
$\sigma_k$

\STATE K-MACE$(x,m_{min}...m_{max},\sigma_{i_{min}}...\sigma_{i_{max}})$ in
feature space (replacing $kmeans(x,m)$ with the kernel
$kmeans(x,m,\sigma)$ clustering algorithm and (\ref{eq:8}) with
(\ref{eq:18000}))

\STATE Use ($\ref{sigma_grad}$) to obtain the optimum $\sigma_i$, $\hat \ sigma$, the final
clustering solution and $\hat{m}$

\end{algorithmic}
\end{algorithm}

\section{Simulations and Results}

\begin{figure}
\centering
\begin{subfigure}[b]{0.23\textwidth}
    \centering
    \includegraphics[width=.8\linewidth]{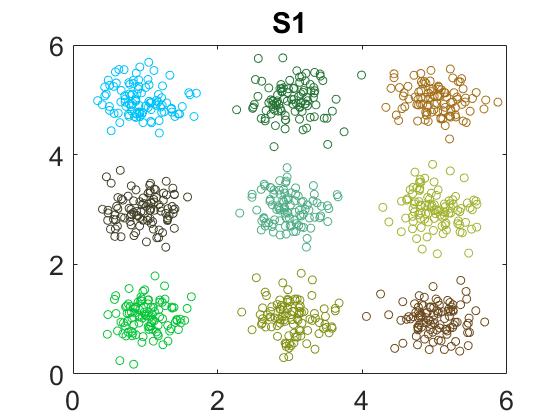} 
    \caption{S1: Uniform variance \\ and no overlaps ($N=900$)} 
    \label{cv1}
\end{subfigure}
\begin{subfigure}[b]{0.23\textwidth}
    \centering
    \includegraphics[width=.8\linewidth]{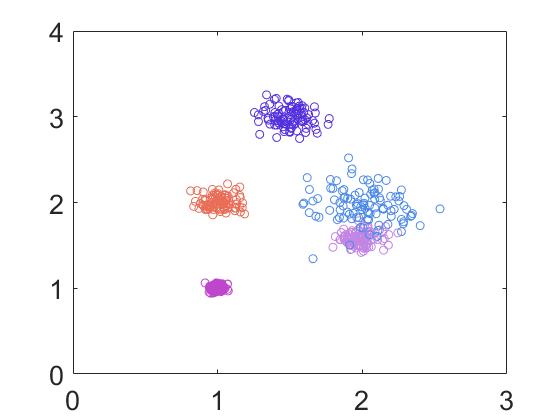}
     \caption{S2:Varying variance \\ and major overlaps  ($N=500$)}
       \label{cv2}
\end{subfigure}
\vskip\baselineskip
\begin{subfigure}[b]{0.23\textwidth}
    \centering
    \includegraphics[width=.8\linewidth]{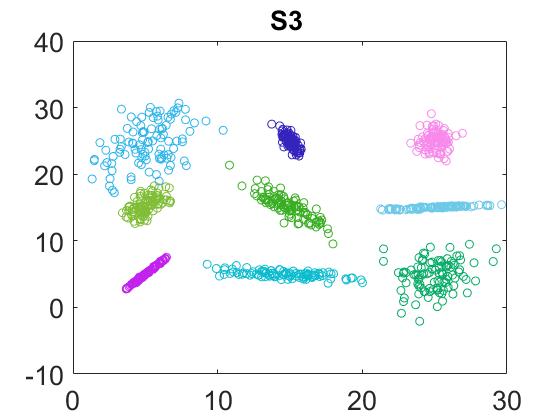}
    \caption{S3: Varying covariance \\ and minor overlaps ($N=900$)} 
     \label{cv3} 
\end{subfigure}
\begin{subfigure}[b]{0.23\textwidth}
    \centering
    \includegraphics[width=.8\linewidth]{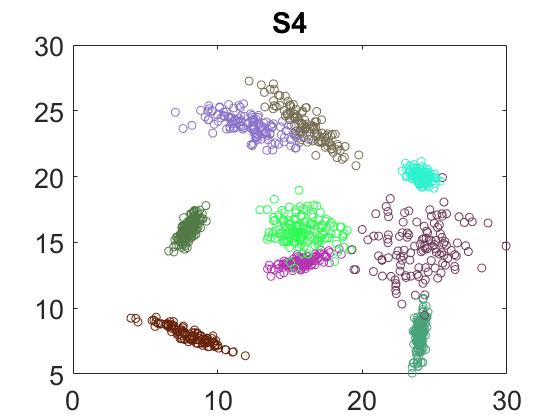}
    \caption{S4: Varying covariance \\ and major overlaps ($N=500$)} 
      \label{cv4}
\end{subfigure}
\caption{Synthetic Data}
\end{figure}

\begin{figure} 
\centering
\includegraphics[scale=0.2]{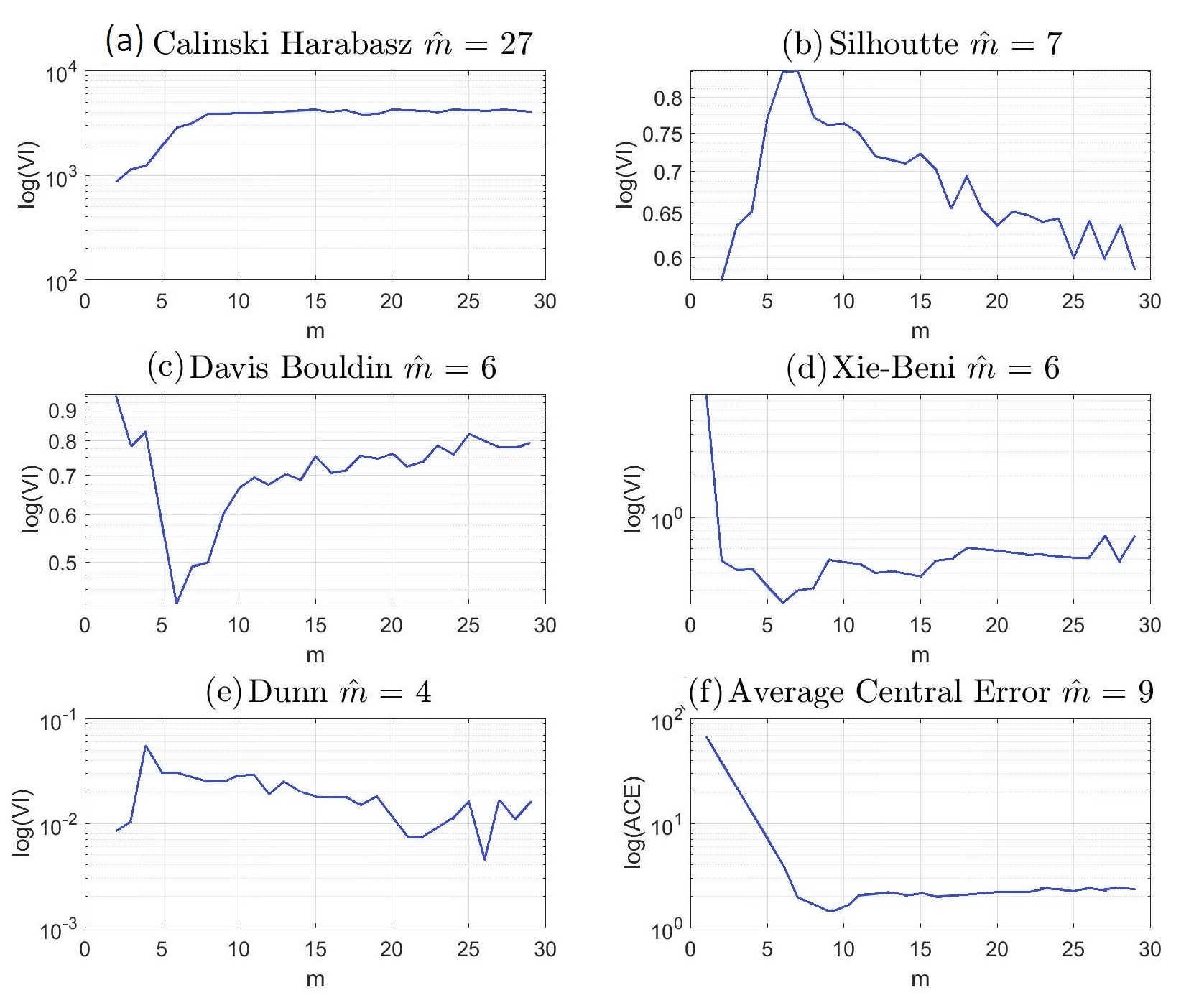}
\caption{ACE and validity index values for data set S4. For each approach Optimum cluster number happens at $m=\hat{m}$. ($\overline{m}=9$)}
\label{figure:0}
\end{figure}
We compare the proposed methods with other well known methods for several synthetic and real data sets. K-MACE is compared with its predecessor MACE-means \cite{ref:9} as well as other well known index validity methods that are CH+K-means \cite{CB}, DB+K-means \cite{DB}, Sil+K-means\cite{sil} and gap+K-means \cite{gap}. We also compare K-MACE and kernel-KMACE (using Gaussian Kernel) with two divisive hierarchical clustering methods that are also partitioning clustering schemes, G-means \cite{gmeans} that is mainly proposed for Gaussian clustering, as well as Dip-means \cite{dip-means} the is a well known most recent approach. Another class of partitioning clustering approach is the density based clustering. We compare the methods with DBSCAN \cite{dbscan} and mean-shift \cite{mean_shift}, two of the most well known algorithm under this category. We also include OPTICS \cite{optics} which is an improved version of the DBSCAN. Note that OPTICS does not provide the clustering solution but only provides the number of clusters estimate. Results of X-means are worse than dip-means and G-means for all our simulation results and therefore, we do not provide them. This observation is consistent with what has already been claimed about G-means and dip-means outperforming X-means \cite{gmeans}, \cite{dip-means}. 

For the purpose of evaluating the performance of clustering algorithms, we record the average and the standard deviation of the value of the CNC estimate, $E[\hat{m}] \pm std[\hat{m}]$, as well as the accuracy of the CNC estimate. The accuracy is the percentage that the algorithm identifying the CNC correctly. The Adjusted Random Index (ARI) \cite{ari} and the Normalized Variation Information (NVI) \cite{nvi} measure agreement between the estimated clustering solution to that of the true partition. ARI scores close to $100\% $ and NVI scores close to $0\%$ indicate full agreement between the algorithm's clustering solution and  the true partition.

\subsection{Synthetic Data}
The first set of generated 2-D data are shown in Figures \ref{cv1} - \ref{cv4}. These data sets have varying complexity in terms of degrees of overlap, type of cluster distribution, and number of elements. To generate these data different samples of scatter factors are generated 100 times around fixed centers $\bar c$s. While S1 has nine non-overlapped clusters, there is a major overlapping of two of the five clusters of S2. clusters in S1 and S2 have uncorrelated (spherical) scatter factors, i.e., $\Sigma(i)$s in (\ref{eq:added_1}), are diagonal matrices. However, S3 and S4 have non-diagonal $\Sigma(i)$s (ellipsoidal) and we denote that as correlated scatter factors for each cluster. Clusters in S3 have some overlap while clusters in S4 have relatively much more overlap. 

Figure \ref{figure:0} shows validity index values of different approaches that are minimum (or maximum) at their estimated CNC.  As the figure shows, the ACE criterion, $z_{m}$, in Figure \ref{figure:0}(f) is the most reliable validity index in this case which points to the correct number of clusters by minimization.

\begin{figure}
\centering
\begin{subfigure}{0.45\textwidth}
	\centering
    \includegraphics[width=1\linewidth]{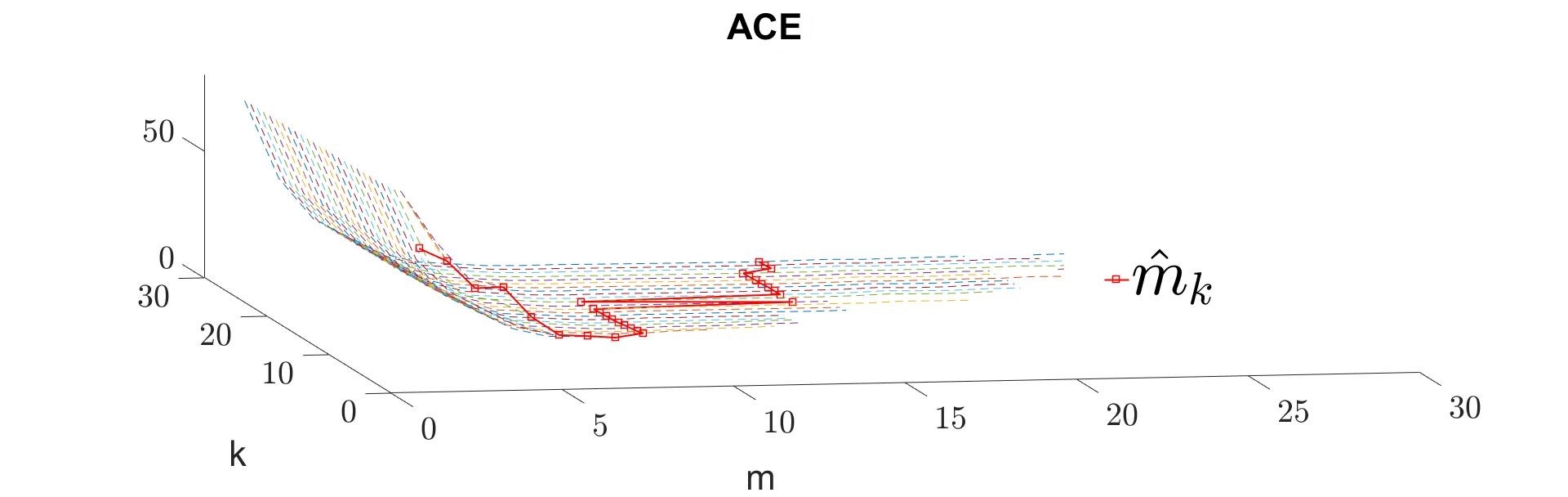} 
    \caption{3D View} 
    \label{3dview}
\end{subfigure}
\begin{subfigure}{0.45\textwidth}
	\centering
    \includegraphics[width=1\linewidth]{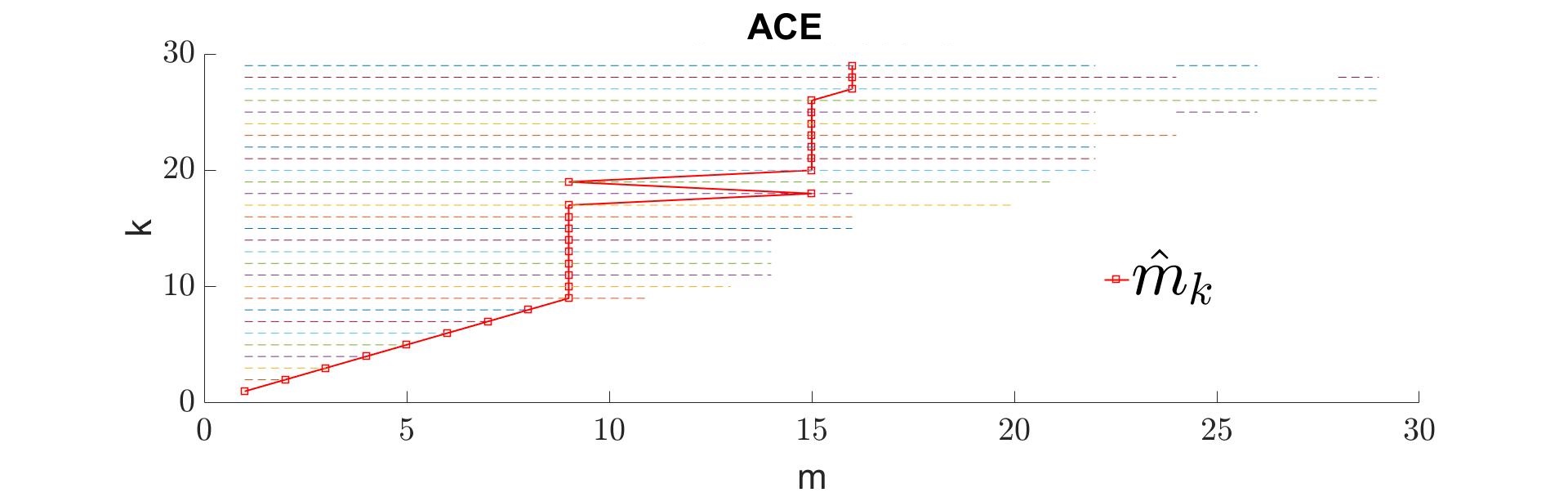} 
    \caption{Top View} 
    \label{topview}
\end{subfigure}
 \caption{Example of $\overline{z_{m,k}}$ behavior as a function of $k$ and $m$ for data set S4 shown in Figure \ref{cv4}} 
 \label{zsm3d}
\end{figure}

Figure \ref{zsm3d} shows how Algorithm \ref{kmace_alg} is implemented through a 3-dimensional plot to visualize behavior of $\overline{z_{m,k}}$ with respect to $k$ and $m$ for S4. The red line shows the values of $\hat {m}_k$ in (\ref{mk}) that minimize ACE for different values of $k$. The optimum value of $k^*$, in (\ref{hhhhh}), in this case is 9 for which its $\hat m$ is also 9. 

\begin{table} 
\centering
\scriptsize
\caption{First Synthetic Data results from average of 100 runs. Results are in the form of $E[\hat{m}] \pm std[\hat{m}]$}
\label{art_result}
\begin{tabular}{|l|l|l|l|l|}
\hline
        &  \begin{tabular}[c]{@{}l@{}}S1\\ m=9\end{tabular} &  \begin{tabular}[c]{@{}l@{}}S2\\ m=5\end{tabular}  &  \begin{tabular}[c]{@{}l@{}}S3\\ m=9\end{tabular}  &  \begin{tabular}[c]{@{}l@{}}S4\\ m=9\end{tabular} \\ \hline
Kernel K-MACE           		  		& {$\mathbf{ 9\pm0}$}		 &		 {$\mathbf{5\pm0}$}		 &		 {$\mathbf{ 9\pm0}$}		 & 		{$\mathbf{ 9\pm0}$}		\\ 
Accuracy (\%)					& {$\mathbf{100}$}			& {$\mathbf{100}$}				& {$\mathbf{100}$}					& {$\mathbf{100}$}			\\
ARI,NVI(\%)						& {$\mathbf{100,0}$}			& {$\mathbf{92,9}$}			& {$\mathbf{98,4}$}					& {$\mathbf{88,9}$}		
\\ \hline
K-MACE           		  		& {$\mathbf{9\pm0}$}		 &		{$\mathbf{5\pm0}$}		 &		 $\mathbf{9\pm0}$		 & 		$\mathbf{9\pm0}$		\\ 
Accuracy (\%)					& {$\mathbf{100}$}			&$\mathbf{100}$				& $\mathbf{100}$					& $\mathbf{100}$			\\
ARI,NVI(\%)						& {$\mathbf{100,0}$}			& $\mathbf{90,10}$			& $\mathbf{95,5}$					& $\mathbf{86,12}$		
\\ \hline

MACE-means       				& {$\mathbf{ 9\pm0}$} 		& 		{$ 7\pm0.8$} 		& 		{$ 12.1\pm0.3$} 			& 		{$ 14\pm0 $}		 	\\ 
Accuracy (\%)					& {$\mathbf{100}$}					& {${0}$}				& {${0}$}								& {${0}$}			\\
ARI,NVI(\%)						& {$\mathbf{100,0}$}				& {${70,30}$}			& {${70,30}$}							& {${80,20}$}	
\\ \hline

Dipmeans           				& {$\mathbf{ 9\pm0}$} 			&		 {$ 4.2\pm0.2 $} 			& {$9\pm0$}					& 		{$ 7\pm0 $}		 \\
Accuracy (\%)					& {$\mathbf{100}$}					& {${65}$}					& {$100$}					& {${0}$}				\\
ARI,NVI(\%)						& {$\mathbf{100,0}$}				& {${76,27}$}				& {$95,5$}					& {${70,30}$}	
\\ \hline       	

Gmeans            				& {$\mathbf{ 9\pm0}$} 			&		 {$ 6\pm0 $} 		& 		{$ 11\pm0 $} 			& 		{$ 17\pm0 $}					 \\ 
Accuracy (\%)					& {$\mathbf{100}$}						& {${0}$} 				& {${0}$}				& {${0}$}							\\
ARI,NVI(\%)						& {$\mathbf{100,0}$}				& {${70,30}$}				& {${80,20}$}			& {${80,20}$}	
\\ \hline

DBSCAN            				& {$\mathbf{ 9\pm0}$} 		&		 {$ 2\pm0 $} 		&		 {$ 11\pm0$} 			&		 {$ 11\pm0    $}			 \\ 
Accuracy (\%)					& {$\mathbf{100}$}					& {${0}$}				& {${0}$}							& {${0}$}				\\
ARI,NVI(\%)						& {$\mathbf{100,0}$}			& {${0,100}$}				& {${80,20}$}						& {${80,20}$}	
\\ \hline

CH + K-means      				& {$\mathbf{ 9\pm0}$}		 &		 {$4.9\pm0.25$} 		&		 {$ 17\pm0$} 			& 		{$ 17.4\pm1.30 $} 			\\ 
Accuracy (\%)						& {$\mathbf{100}$}				& {${85}$}						& {${0}$}						& {${0}$}							\\
ARI,NVI(\%)						& {$\mathbf{100,0}$}				& {${80,20}$}					& {${70,20}$}					& {${70,30}$}		
\\ \hline

Sil + K-means     				& {$\mathbf{ 9\pm0}$}		&		 {$ 4\pm0$} 			& 		{$9\pm0$} 		& 		{$7.8\pm0.4   $}			 \\
Accuracy (\%)					& {$\mathbf{100}$}				& {${0}$}								& {$100$}			& {${0}$}				\\
ARI,NVI(\%)						& {$\mathbf{100,0}$}			& {${80,20}$}						& {$95,5$}				& {${80,20}$}	
\\ \hline

DB + K-means      				& {$\mathbf{ 9\pm0}$}		 &		 {$ 4\pm0$} 			&		 {$9\pm0$}		 &  		{$ 7.7\pm0.5   $} 			\\
Accuracy (\%)						& {$\mathbf{100}$}					& {${0}$}						& {$100$}				& {${0}$}					\\
ARI,NVI(\%)						& {$\mathbf{100,0}$}					& {${80,20}$}				& {$95,5$}					& {${80,20}$}	
\\ \hline

gap + K-means      				& {$\mathbf{ 9\pm0}$}		 &		 {$5\pm0$} 			&		 {${ 14.3\pm2.3}$}		 &  		{$15.9\pm1.4$} 			  \\
Accuracy (\%)						& {$\mathbf{100}$}					& {$100$}				& {${0}$}						& {${0}$}					\\
ARI,NVI(\%)						& {$\mathbf{100,0}$}				& {$90,10$}					& {${80,20}$}					& {${70,20}$}
\\ \hline

Mean-Shift     				& {$\mathbf{ 9\pm0}$}		&		 {$ 4\pm0$} 			& 		{$9\pm0$} 		& 		{$7.0\pm0.2   $}			 \\
Accuracy (\%)					& {$\mathbf{100}$}				& {${0}$}							& {${100}$}			& {${0}$}					\\
ARI,NVI(\%)						& {$\mathbf{100,0}$}			& {${80,20}$}				& {$95,5$}						& {${70,20}$}
\\ \hline

Optics      				& {$\mathbf{9\pm0}$}		 &		 {$ 7\pm0$} 			&		 {$11\pm0$}		 &  		{$ 8\pm0   $} 			\\
Accuracy (\%)				& {$\mathbf{100}$}				& {${0}$}						& {${0}$}				& {${0}$}		\\
ARI,NVI(\%)						&NA							&NA							&NA					&NA
\\ \hline
\end{tabular}
\end{table}

Table \ref{art_result} provides comparison of result for S1-S4 data sets. While all methods provide perfect results for S1 that was a well behaved set of clusters, only few methods are comparable with K-MACE and Kernel K-MACE for S2 and S3 and for data set S4, that is the most complicated set with relatively more overlapping, the only two methods that outperforms other methods are the K-MACE ones. Note that for all these scenarios both  K-MACE approaches standard deviation are zero which indicates precision and robustness of the algorithm to change of scatter factor samples. Not only the K-MACE methods are successful in CNC estimation but also they are providing the highest value of ARI and lowest of NVI that indicates better results of clustering as well. 

The second set of generated data are six clustering sets of data U1-U6 of length 900 with dimension ten that is a higher dimension compared to the first set set S1-S4. The data set is generated by randomly choosing 9 centers in a hypercube for each of the 100 runs. For all these data sets the standard deviation for scatter factors per dimension is smaller than 0.3. For U1, U3 and U5 the hypercube side length is 10 while for U2, U4 an U6 the hypercube side length is shrunk to smaller value of 8. Consequently, the latter sets of data have much greater chance of more overlapping than the first sets.

Simulation results are in Table \ref{art_result2}.  
As the table shows, while for the easier case of U1 almost all methods are successful, for the other data sets the only methods that consistently outperform other methods are the K-MACE ones. Mean-shift and Dip-means were able to estimate CNC correctly for some  data sets . It is work mentioning that Kernel k-MACE shows its expected superiority for the cases with major   overlap. While K-MACE is compatible with kernel K-MACE for U2, for U4 and U6 Kernel K-MACE performs better and with a very high accuracy.  
 
\begin{figure}
\centering
\begin{subfigure}{0.4\textwidth}
	\centering
    \includegraphics[width=1\linewidth]{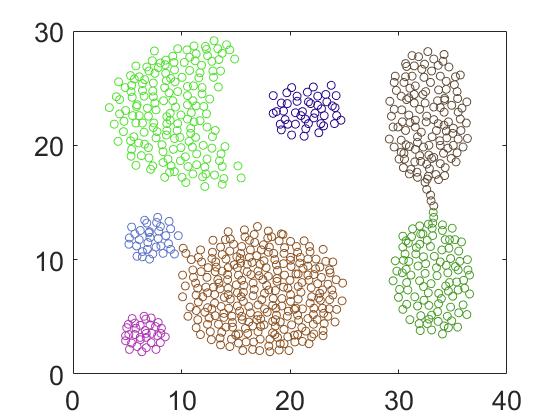} 
    \caption{Aggregation Data Set} 
    \label{aggregation}
\end{subfigure}
\begin{subfigure}{0.4\textwidth}
	\centering
    \includegraphics[width=1\linewidth]{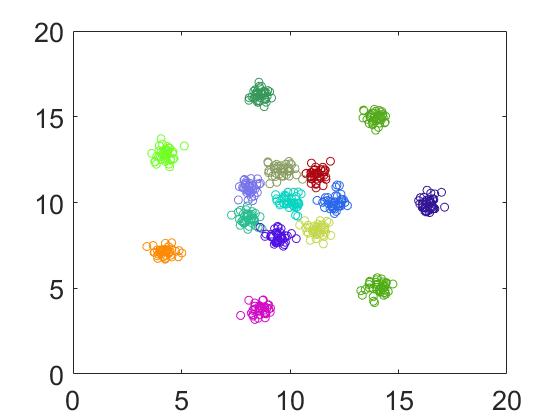} 
    \caption{S15 Data Set} 
    \label{s15}
\end{subfigure}
\begin{subfigure}{0.4\textwidth}
	\centering
    \includegraphics[width=1\linewidth]{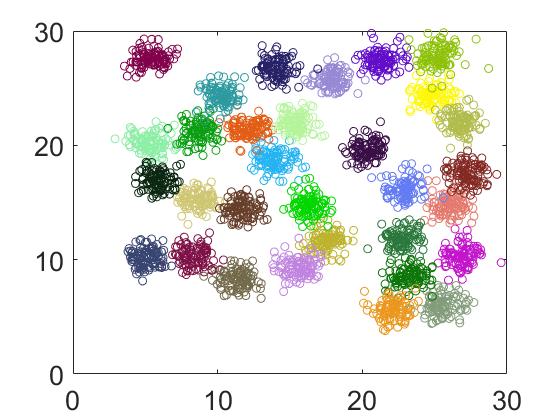} 
    \caption{D31 Data Set} 
    \label{d31}
\end{subfigure}
 \caption{Synthetic Data Set. } 
\end{figure}

Third set of synthetic data are from \cite{d31}- \cite{s15} and depicted in Figures \ref{aggregation} - \ref{d31}. As the figures show while Aggregation data set has non Gaussian shaped clusters, clusters of  S15 have very few points and 31 clusters of D31 are touching.  Simulation results  for this data are presented in Table \ref{art_result3}. As the table shows, while K-MACE approaches choose 6 clusters for Aggregation instead of 7, nevertheless are still the winners among all the methods as still are the closest CNC estimate and more importantly provide  large ARIs and low NVIs. The next method that competes with K-MACE methods for this data is DBSCAN. As S15 is a much easier data set with a clean clustering structure, couple of methods, including K-MACE, perform equally well for this data set. For D31, K-MACE and especially Kernel K-MACE are on the average outperforming other methods even though Gmeans is the only other method that provides the correct CNC estimate. This example shows importance of overall performance of the clustering method. As the table shows, kernel K-MACE has $95\%$ accuracy which is the largest accuracy among the methods after Gmeans. In addition, it has the highest AVI and the lowest NVI. It is worth mentioning that all clustering methods have tuning parameters that we set to their default values \cite{dip-means}. For example, DBSCAN has two parameters, $\epsilon$ and $MinPts$ that indirectly control the number of clusters. it is known that DBSCAN is fairly sensitive to these parameters and if not finely tuned, can result into bad clustering solution. 
Another example is Gmean that is using `Anderson-Darling threshold' \cite{gmeans} for the splitting procedure. This value is by default $0.95$. In Table \ref{art_result3} results of changing this value to $0.7$ is also shown. As the table shows the Gmeans performance is degraded with this new value. K-MACE approaches have parameters $\alpha_{mj}$ and $\beta$. The default value that we used is five for all these parameters. This value guarantees confidence and validation probabilities to be larger than $0.96=1-\frac{1}{25}$. Note that K-MACE is not sensitive to changing the default values. Changing this value to larger numbers, as long as the condition in Section \ref{vall} is satisfied, still provides the same results of this default value with a very small variation. 

\begin{table*} 
\centering
\scriptsize
\caption{Second Synthetic Data Set, nine cluster sets ($\overline{m}=9$) with higher dimension ($d=10$) and $N=900$, averaged over 100 runs. Results are in the form of $E[\hat{m}] \pm std[\hat{m}]$. (Uncorr-Scatt stands for uncorrelated (spherical) scatter factors, Corr-Scatt stands for Correlated (ellipsoidal) scatter factors)} 
\label{art_result2}
\begin{tabular}{|l|l|l|l|l|l|l|}
\hline
          &  \begin{tabular}[c]{@{}l@{}}\textbf{U1}\\ Uncorr-Scatt \\ Same variance \\minor overlap \end{tabular} &  \begin{tabular}[c]{@{}l@{}}\textbf{U2}\\ Uncorr-Scatt\\ Same variance \\major overlap\end{tabular}  &  \begin{tabular}[c]{@{}l@{}}\textbf{U3}\\Uncorr-Scatt\\ Different variance \\minor overlap\end{tabular}  &  \begin{tabular}[c]{@{}l@{}}\textbf{U4}\\Uncorr-Scatt\\ Different variance \\ major overlap\end{tabular}  &  \begin{tabular}[c]{@{}l@{}}\textbf{U5}\\Corr-Scatt\\ Different Covariance \\ minor overlap\end{tabular}  &  \begin{tabular}[c]{@{}l@{}}\textbf{U6}\\Corr-Scatt\\ Different Covariance \\ major overlap\end{tabular}  \\ \hline
Kernel K-MACE           		  		& {$\mathbf{ 9\pm0}$}		 &		 {$\mathbf{9\pm0}$}		 &		 {$\mathbf{ 9\pm0}$}		 & 		{$\mathbf{ 9\pm0}$}	 & 		{$\mathbf{9\pm0}$} & 	{$\mathbf{ 9\pm0}$}		\\ 
Accuracy (\%)				& {$\mathbf{100}$}					& {$\mathbf{100}$}					& {$\mathbf{100}$}					& {$\mathbf{95}$}			& {$\mathbf{100}$}			& {$\mathbf{96}$} 		\\

ARI,NVI(\%) & {$\mathbf{98,2}$} & {$\mathbf{93,8}$} & {$\mathbf{96,5}$} &
{$\mathbf{95,5}$} & {$\mathbf{96,3}$} & {$\mathbf{95,5}$}
\\ \hline

K-MACE           		  		& {$\mathbf{ 9\pm0}$}		 &		 {$\mathbf{9\pm0}$}		 &		 {$\mathbf{9\pm0}$}		 & 		{$ 8.6\pm0.7$}	 & 		{$\mathbf{9\pm0}$} & 	{$8.8\pm0.1$}		\\ 
Accuracy (\%)				& {$\mathbf{100}$}					& {$\mathbf{100}$}					& {$\mathbf{100}$}					& {$92$}			& {$\mathbf{100}$}			& {$95$} 		\\
ARI,NVI(\%)					& {$\mathbf{95,5}$}			& {$\mathbf{90,10}$}				& {$\mathbf{90,5}$}					& {$90,10$}			& {$\mathbf{95,5}$}		& {$95,10$}
\\ \hline

MACE-means       				&{${ 11\pm0}$}		& 		{$14.1\pm0.2$} 		& 		 {${ 10\pm0.2}$}			& 		{$12.8\pm1.3$}  & 				{$11.8\pm1.2$}			& {$13.5\pm2.1$}	 	\\ 
Accuracy (\%)					& {$0$}					& {${0}$}				& {${5}$}								& {${12}$}							& {${15}$}			& {${8}$}	\\
ARI,NVI(\%)						& {$90,10$}				& {$80,20$}					& {${90,10}$}						& {$80,10$}							& {$70,20$}			& {$80,20$}	
\\ \hline

Dipmeans           				& {$ 9\pm0$} 			&		{$9\pm0$} 		& 		 {${  8.6\pm1.4}$}					& 		{$ 6.84\pm2.62 $}	 & 		{$9\pm0$}	  & {$  8.56\pm1.51 $}	 \\ 
Accuracy (\%)					& {$100$}						& {$100$}			& {${69}$}							& {${61}$}							& {$100$}			&	  {${76}$}	\\
ARI,NVI(\%)						& {$95,5$}					& {$90,10$}			& {$80,20$}						& {$70,20$}						& {$90,10$}			& {$80,20$}	
\\ \hline

Gmeans            				& {$9\pm0$} 			&		 {$ 9.16\pm0.42 $} 		& 		{$ 9.12\pm0.38 $} 			& 		{$ 12.32\pm4.37 $}		& 		{$ 9.12\pm0.38 $}		& 		{$ 9.14\pm0.35 $}			 \\ 
Accuracy (\%)					& {$100$}			&		  {${88}$} 				&		  {${92}$} 					&		 {${17}$} 				&		 {${88}$}				&		 {${91}$}\\
ARI,NVI(\%)						& {$95,5$}				& {${90,10}$}				&{90,10}					& {$80,10$}							&{90,10}					&	 {$90,10$}	
\\ \hline

DBSCAN            				&  {$8.8\pm0.4$} 		&		 {$8.9\pm0.4$} 		&		 {$6.9\pm1.6$} 			&		 {$ 11\pm0    $}		&  		{$8.9\pm0.2$}		& 		{$8.7\pm0.5$}	 \\
Accuracy (\%)				& {${87}$}		&        			 {${96}$}		& 				{${0}$}		& 				{${0}$}			& 				{${95}$}	& 				{${79}$}\\
ARI,NVI(\%)					& {${90,10}$}					&{90,10}						&{70,25}					& {$90,10$}					&{90,10}				&{90,10}	
\\ \hline

CH + K-means      				& {$9\pm0$}		 &		{$9.1\pm0.3$} 		&		 {${9.1\pm0.3}$} 			& 		{$ 17.4\pm1.30 $} 	 & 		{${ 18\pm1.2}$}	 & 		{${ 18.7\pm0.8}$}		\\
Accuracy (\%)					& {$100$}			& 		{${90}$}			& 			{${90}$}						& 		{${0}$}				& {${0}$}	 & 			{${0}$}	\\
ARI,NVI(\%)						& {$95,5$}					&{95,5}						&{90,10}							& 	{$80,30$}				& {${40,50}$}			&{40,50}	
\\ \hline

Sil + K-means     				& {$9\pm0$}		&		 {$9\pm0.3$}			& 		{$ 9\pm0.3$} 			& 		{$7.8\pm0.4   $}		& 		{$ 9\pm0$} 	 & 		{${ 8.2\pm1.1}$}	  \\
Accuracy (\%)					& {$100$}		&   		 {${96}$}					& 			{${96}$}			&		 {${0}$}			& 		{$100$}		 & 		{${68}$}\\
ARI,NVI(\%)						& {$95,5$}				&{95,5}								&{95,2}					& 	{$80,20$}							& {$90,10$}		&{90,10}	
\\ \hline

DB + K-means      				& {$9\pm0$}		 &		 {$8.9\pm0.3$}	 		&		  {$ 8.9\pm0.3$}		 &  		{$ 7.7\pm0.5   $} 	& 		{$9\pm0$}	    & 		{${ 7.1\pm1.5}$}		\\
Accuracy (\%)					& {$100$}			&   		 {${91}$}				&   		 {${91}$}			&   		 {${6}$}			&   		 {$100$}				&   		 {${2}$}\\
ARI,NVI(\%)						& {$95,5$}				&{95,5}								&{90,10}					& 	{$80,20$}				& {$90,10$}			&{80,10}	
\\ \hline

gap + K-means      				& {$9\pm0$}		 &		{$9.1\pm0.3$} 			&		 {${9.1\pm0.3}$} 		 &  		{$15.9\pm1.4$} 	& 		{$11.8\pm5.5$} 	 & 		{${ 13.1\pm0.8}$}		  \\ 
Accuracy (\%)						& {$100$}		&   		 {${91}$}				&   		 {${91}$}				&   		 {${0}$}		&   		 {${6}$}			&   		 {${0}$}\\
ARI,NVI(\%)						& {$95,5$}					&{95,7}							&{90,10}					& 	{$70,30$}				& 	{$80,30$}				&{80,25}	
\\ \hline

Mean-Shift     				& {$9\pm0$}		&		 {$ 9\pm0$} 			& 		 {${ 9\pm0.1}$} 		& 		{$7.0\pm0.2   $}	 	& 	{$9\pm0$}		&	 {$8.9\pm0.2$}	 \\ 
Accuracy (\%)					& {$100$}					& {$100$}		&   		 {${96}$}					&   		 {${18}$}			& {$100$}		&	 {${96}$}\\
ARI,NVI(\%)						& {$95,5$}					& {$90,10$}				&{95,5}					& 	{$80,20$}				& {$90,10$}		&{90,10}	
\\ \hline

Optics      				& {$9.0\pm1.8$}		 &		 {$9.0\pm1.5$} 			&		 {$9.1\pm1.4$}		 &  		{$ 8\pm0   $} 		& 	{$9.6\pm1.3$}		& 		{$9.4\pm1.9$}	\\ 
Accuracy (\%)			& {${80}$}	& 		{${79}$}				&		 {${82}$}			&		 {${0}$}			&	{${83}$}			& 		{$74$}	\\
ARI,NVI(\%)				&NA						&NA							&NA							&NA					&NA						&NA
\\ \hline

\end{tabular}
\end{table*}

\begin{table}
\centering
\scriptsize
\caption{Third set of Synthetic data results from averaging of 100 runs. Results are in the form of $E[\hat{m}] \pm std[\hat{m}]$}
\label{art_result3}
\begin{tabular}{ | c | c | c | c |}
\hline
{ \ \ \ \ \ \ \ \ \ \ \ \ \ \ \  }  &{aggregation} & {s15} & {D31} \\ 
{ \ \ \ \ \ \ \ \ \ \ \ \ \ \ \  }  &{$\overline{m}=7$} & {$\overline{m}=15$} & {$\overline{m}=31$} \\ \hline

Kernel K-MACE    				& {$\mathbf{6\pm0}$} 		& 		{$ \mathbf{15\pm0}$} 		& 		{$\mathbf{31\pm0}$} 			\\ 
Accuracy (\%)			& {$\mathbf{0}$}					& 		{$\mathbf{100}$}			& 		{$\mathbf{95}$}					\\
ARI,NVI(\%)				& {$\mathbf{91,12}$}		& {$\mathbf{100,0}$}					& {$\mathbf{94,7}$}		
\\ \hline

K-MACE    				& {$\mathbf{6\pm0}$} 		& 		{$ \mathbf{15\pm0}$} 		& 		{${30.9\pm1.9}$} 			\\ 
Accuracy (\%)			& {${0}$}					& 		{$\mathbf{100}$}			& 		{${88}$}					\\
ARI,NVI(\%)				& {$\mathbf{90,10}$}		& {$\mathbf{100,0}$}					& {${90,10}$}		
\\ \hline

MACE-means       		& {${ 16\pm0.6}$} 			& 		{$\mathbf{15\pm0}$} 		& 		{$ 31.8\pm0.6$} 			\\  
Accuracy (\%)			& {${0}$}					& 		{$\mathbf{100}$}			& 		{${72}$}				\\
ARI,NVI(\%)				& {${40,40}$}				& {$\mathbf{100,0}$}					& {${90,10}$}		
\\ \hline

Dipmeans           			& {$5\pm0$} 				&		 {$ 8\pm0 $} 		& 		{$ 23\pm0 $} 						 \\  
Accuracy (\%)			& {${0}$}					& 		{$0$}				& 		{${0}$}				\\
ARI,NVI(\%)				& {${80,20}$}				& 		{${70,30}$}					& {${86,12}$}		
\\ \hline

Gmeans($\alpha = 0.95$) & {$ 9\pm0$} 				&		 {${16\pm0}$} 		& 		{$ {31\pm0}$}			 \\ 
Accuracy (\%)			& {${0}$}					& 		{${0}$}				& 		{${100}$}				\\
ARI,NVI(\%)				& {${60,30}$}				& 		{${95,6}$}					&{${86,12}$}		
\\ \hline

Gmeans($\alpha = 0.7$)     & {$ 13\pm0$} 				&		 {$19.9\pm0.35 $} 		& 		{${ 46\pm2.9}$}			 \\ 
Accuracy (\%)			& {${0}$}					& 		{$0$}				& 		{${0}$}				\\
ARI,NVI(\%)				& {${60,30}$}				& 		{${80,20}$}					& {${86,12}$}		
\\ \hline

DBSCAN            			& {${5.2\pm0.2}$} 			&		 {$ \mathbf{15\pm0} $} 		&		 {$ 34\pm0$} 					 \\  
Accuracy (\%)			& {${0}$}					& 		{$\mathbf{100}$}				& 		{$0$}				\\
ARI,NVI(\%)				& {${90,10}$}				& {$\mathbf{100,0}$}						& {${50,30}$}		
\\ \hline

CH + K-means      		& {${ 15.7\pm0.9}$}		 	&		 {$ 15.3\pm0.5$} 		&		 {$ 3.4\pm1.1$} 					\\  
Accuracy (\%)			& {${0}$}					& 		{$90$}					& 		{$0$}				\\
ARI,NVI(\%)				& {${40,40}$}				&	 {${95,5}$}					& {${90,10}$}			
\\ \hline

Sil + K-means     			& {${ 3\pm0}$}				&		 {$ 15.3\pm0.5$} 			& 		{${30.9\pm1.5}$} 					 \\  
Accuracy (\%)			& {${0}$}					& 		{$90$}						& 		{${83}$}				\\
ARI,NVI(\%)				& {${70,40}$}				&	 {${95,5}$}					& {${90,10}$}		
\\ \hline

DB + K-means      		& {${ 4\pm0}$}				&		 {$ 12.2\pm3.5$} 			&		 {${ 28.4\pm2.0}$}				\\  
Accuracy (\%)			& {${0}$}					& 		{$35$}						& 		{$0$}				\\
ARI,NVI(\%)				& {${80,30}$}				& 		{${90,10}$}					& {${90,10}$}		
\\ \hline

gap + K-means      		& {$14.5\pm0.8$}			&		 {$ 12\pm2.5$} 			& {${ 33.5\pm1.3}$}				\\  
Accuracy (\%)			& {${0}$}					& 		{$18$}					& {$60$}				\\
ARI,NVI(\%)				& {${40,40}$}				&	      {${90,10}$}				& {${90,10}$}		
\\ \hline

Mean-shift      			& {${5.5\pm0.8}$} 		 	&{$ \mathbf{15\pm0}$}  		&		 {${ 12\pm1.8}$}				\\  
Accuracy (\%)			& {${0}$}					& {$\mathbf{100}$}			& 		 {$0$}				\\
ARI,NVI(\%)				& {${70,40}$}				& {$\mathbf{100,0}$}						& {${20,70}$}		
\\ \hline

Optics      				& {${ 8\pm0}$}		 		&		 {$ {15\pm0}$} 		&		 {${ 16.8\pm2.3}$}				\\  
Accuracy (\%)			& {${0}$}					& 		{${100}$}			& 		{$0$}				\\
ARI,NVI(\%)				&NA						&		NA							& 		NA		
\\ \hline		

\end{tabular}
\end{table}

\subsection{Real Data}
Results on eight real-world data sets from UCI machine learning repository website \footnote{http://archive.ics.uci.edu/ml/}, that are seeds, iris, vertebral, wine, breast, Wisconsin Diagnostic Breast Cancer, Ecoli and Multiple Features (dutch handwritten), are shown in Table \ref{real1_result}. 
Values of $\overline{m}$, the actual number of clusters, $d$, the dimension of the data, and $N$, the number of samples in the data set are provided in the table. As the table shows K-MACE methods outperform most of methods. For all data sets K-MACE methods provide accurate estimate of CNC. Even for Ecoli with 8 clusters, the best method among all the methods are the K-MACE methods. However, the CNC estimate in this case is close to 5 by K-MACE. The reason is the structure of the data for which the number of elements, $n_{mj}$, of each cluster are [143, 77, 52, 35, 20, 5, 2, 2]. As these numbers show, the last three clusters have very few elements that can't be detected as independent clusters by K-MACE methods. This is perhaps similarly effecting the other methods such as G-means.

\begin{table*}
\scriptsize
\centering
\caption{Real Data examples (average of 50 runs)}
\label{real1_result}
\newlength\q
\setlength\q{\dimexpr .05\textwidth -1 \tabcolsep}
\noindent
\begin{tabular}{ | c | c | c | c | c | c | c | c | c |}
\hline
{ \ \ \ \ \ \ \ \ \ \ \ \ \ \ \  } &{\begin{tabular}[c]{@{}l@{}} Seeds\\$\overline{m}$=3, d=7\\N=210\end{tabular}}  & {\begin{tabular}[c]{@{}l@{}} Iris\\ $\overline{m}$=3, d=4\\ N=150\end{tabular}}  & {\begin{tabular}[c]{@{}l@{}} Vertebral \\ $\overline{m}$=3,d=6\\ N=310\end{tabular}} & {\begin{tabular}[c]{@{}l@{}} Wine\\ $\overline{m}$=3, d=13\\ N=178\end{tabular}} & {\begin{tabular}[c]{@{}l@{}} Breast\\ $\overline{m}$=2,d=9 \\ N=698\end{tabular}}  & {\begin{tabular}[c]{@{}l@{}} WDBC\\ $\overline{m}$=2,d=32 \\ N=569\end{tabular}} & {\begin{tabular}[c]{@{}l@{}} Ecoli\\ $\overline{m}$=8,d=7 \\ N=335\end{tabular}} & {\begin{tabular}[c]{@{}l@{}} MF\\ $\overline{m}$=10,d=216 \\ N=2000\end{tabular}}  \\ \hline

{Kernel K-MACE}        				& {$\mathbf{3\pm0}$}   & {$\mathbf{3\pm0}$}&{$\mathbf{3\pm0}$} &{$\mathbf{3\pm0}$} &  {$\mathbf{2\pm0}$}  &{$\mathbf{2\pm0}$} &{$\mathbf{5\pm0}$} &{$\mathbf{10\pm0}$}\\
\hline

{K-MACE}        				& {$\mathbf{3\pm0}$}   & {$\mathbf{3\pm0}$}&{$\mathbf{3\pm0}$} &{$\mathbf{3\pm0}$} &  {$\mathbf{2\pm0}$}  &{$\mathbf{2\pm0}$} &{$4.9\pm0.3$} &{$\mathbf{10\pm0}$}\\      
\hline				
				 				
{MACE-means}    		   		& {$\mathbf{3\pm0}$}		&{$6.5\pm0.5$}           		&{$1\pm0$}           			  &{$11.2\pm0$}		&       {$1\pm0$}       &{$7.2\pm2.1$} &{$1\pm0$} &{$5\pm0$}       \\ 	
\hline

{DIP-means}       				    		 & {$1\pm0$}          			&{$2\pm0$}              		&{$1\pm0$}              		 &{$1\pm0$}   & 	{$12\pm0$}    &{$1\pm0$}    &{$2\pm0$} &{$4\pm0$}\\  	
\hline

{G-means}       				    		 & {$4\pm0$}          			&{$4\pm0$}		             	&{$5\pm0$}              		 &{$2\pm0$}    &	{$96\pm0$}  &{$14\pm0$} &{${4\pm0}$}   &{$35\pm0$}\\  
\hline

{DBSCAN}       				     		 & {$2\pm0$}          			&{$\mathbf{3\pm0}$}             	&{$1\pm0$}              		 &{$1.0\pm0.0$}    &          {$11\pm0$}  &		{$1\pm0$} &       {$1\pm0$}  &       {$1\pm0$}  \\
\hline

{CH + K-means}   			     		& {$\mathbf{3\pm0}$}    	  	&{$\mathbf{3\pm0}$}             	&{$4\pm0$}               		 &{$12.7\pm0.5$}  & {$\mathbf{2\pm0}$}      &{$11.1\pm1.0$} &{$3.0\pm0.0$}&{$4\pm0$}\\    	
\hline
{Sil + K-means}   			 		& {$2\pm0$}            		&{$2\pm0$}            			&{$\mathbf{3\pm0.2}$}          		&{$2.0\pm0.0$}     &    {$\mathbf{2\pm0}$}   &{$2.8\pm0.4$} &{$3.0\pm0.0$}&{$4.9\pm0.3$}\\	
\hline				 				
{DB + K-means}    		  		 &{$2\pm0$}                 		&{$2\pm0$}                   		&{$3.1\pm0.3$}            		  &{$7\pm0.5$}     &           {$\mathbf{2\pm0}$}  & {$\mathbf{2\pm0}$} &{$2.0\pm0.0$}& {$5.0\pm0.0$}	\\	
\hline

{gap + K-means}  				& {$2\pm0$}            		&{$2\pm0$}            			&{$5.2\pm0.5$}            		 &{$8.7\pm5.2$}      &{$10.42\pm0.7$}  & {$10.5\pm1.3$} & {$16.0\pm1.4$} & {$7.0\pm0.0$}	\\
\hline
{mean shift}  			 		& {$3.9\pm0.5$}            		&{$2\pm0$}            			&{$18.7\pm0.9$}            		 &{${3.7\pm0.4}$}   & {$152.5\pm2.2$}  &{$7.6\pm0.6$} &{$22.4\pm0.9$}   & {$112.5\pm2.6$}  \\	
\hline
{optics}  			 	& {$8\pm0.0$}           		& {$5\pm0.0$}           			&{$1.0\pm0.0$}            		 &{$10\pm0.0$}    	&{$9\pm0$}     &{$4.3\pm0.7$} & {$1.0\pm0.0$} & {$1.0\pm0.0$}  \\	
 \hline

\end{tabular}
\end{table*}
\section{Conclusion}
K-MACE algorithm is a validity index method for CNC estimation in K-means clustering. The unique method of ACE estimation uses only the available data and the available cluster compactness. The approach estimates the required scatter factor variance of clusters by evaluating consistency of the $m$-clustering with the prior assumption on the cluster structure. Due to this insightful tactic, the method shows robustness in the clustering procedure. Kernel K-MACE algorithm is an implementation of K-MACE in feature space using kernel functions which produces better clustering results for data sets with a large degree of overlap between clusters. In addition, while existing methods tune the Gaussian kernel parameters by trial and error, the proposed method automatically tunes these parameters. The simulation results also confirm that  K-MACE and kernel K-MACE not only outperform the competing methods in CNC estimation, but also perform with more precision in clustering the data in the sense of ARI and NVI even for clusters with higher level of overlapping. 
\appendices

\section{Proof of Lemma 1}
\label{app_lemma1}

ACE in each cluster in (\ref{zsmj}) can also be written in the following form  
\begin{equation}
\begin{split} 
z_{{mj}}  = \norm{\overline{c}_{x_{mj}} - {x_{mj}} B_{mj}}_2^2  \label{zzz}
\end{split}
\end{equation}
where $x_{mj}$ is a vector of elements of $C_{mj}$, $\bar c_{x_{mj}}$ is a vector of all the associated $\bar c_{x_{mj(i)}}$s to the elements of $C_{mj}$ and 
$B_{mj}$ is an $n_{mj}$ by $n_{mj}$ averaging matrix 
\begin{equation} \label{BBB}
B_{mj} = 
\begin{bmatrix}
\frac{1}{n_{mj}}& \dots{} & \frac{1}{n_{mj}} \\
\vdots &\ddots & \vdots\\
\frac{1}{n_{mj}} &\dots{} & \frac{1}{n_{mj}}
\end{bmatrix}
\end{equation}

Replacing $x_{mj}(i)$ from (\ref{eq:dd_mod}), ($x_{mj}(i) = \overline{c}_{x_{mj}(i)} +\overline{W}_{x_{mj}(i)}$ in (\ref{zzz}) we have: 
\begin{equation}
\begin{split}
Z_{{mj}}  = \norm{  \overline{c}_{x_{mj}}(I-B_{mj})- \overline{W}_{x_{mj}}B_{mj} }_2^2 \\
Z_{{mj}}  = \norm{  \overline{c}_{x_{mj}}A_{mj} - \overline{W}_{x_{mj}}B_{mj} }_2^2 
\end{split}
\end{equation}

and 
$A_{mj}=I-B_{mj}$ is shown in (\ref{eq:11}).
Since $A_{mj}^TB_{mj} = 0$, we have

\begin{equation}\label{app:zsm1}
Z_{{mj}}  = \norm{\Delta_{{mj}}}_2^2 +  \norm{ \overline{W}_{x_{mj}}B_{mj}}_2^2
\end{equation}
where $ \norm{\Delta_{{mj}}}_2^2 = \norm{\overline{c}_{x_{mj}}A_{mj}}_2^2$
\begin{equation}\label{app:zsm2}
\scriptsize
\norm{ \overline{W}_{x_{mj}}B_{mj} }_2^2  =  \frac{1}{n_{mj}}\sum_{i=1}^{n_{mj}}\overline{W}_{x_{mj}(i)}^T \overline{W}_{x_{mj}(i)} + \frac{1}{n_{mj}}\sum_{i \neq k}^{n_{mj}} \overline{W}^T_{x_{mj}(i)} \overline{W}_{x_{mj}(k)}
\end{equation}

Plugging (\ref{app:zsm2}) in (\ref{app:zsm1}) will yield the expression in (\ref{eq:10}). 
The first term of  (\ref{eq:10}) is a constant. The second term is chi-square random variable with a non-zero mean. Because of the independence between the random vectors $\overline{W}_{x_{mj}(i)}$ for $i=[1,2,\dots,N]$, the expectation can be brought inside the summation term and 

\begin{equation}\label{exp_inside}
\begin{split}
E\left[ \frac{1}{n_{mj}}\sum_{i=1}^{n_{mj}}\overline{W}_{x_{mj}(i)}^T \overline{W}_{x_{mj}(i)} \right] =\\
 \frac{1}{n_{mj}}\sum_{i=1}^{n_{mj}}E[\overline{W}_{x_{mj}(i)}^T \overline{W}_{x_{mj}(i)}] = \frac{1}{n_{mj}}\sum_{i=1}^{n_{mj}}tr(\overline{\wedge}_{x_{mj}(i)})
\end{split}
\end{equation}

The third term has an expected value of 0 due to independence between $\overline{W}_{x_{mj}(i)}$s. Thus, the overall expected value of $Z_{{mj}}$ is given by (\ref{eq:12}). Note that the first term is a constant thus, it has a zero variance. The variance of the other two terms are given below:

\begin{equation}\label{var_zsm_app1}
\scriptsize
Var \left[  \frac{1}{n_{mj}}\sum_{i=1}^{n_{mj}}\overline{W}_{x_{mj}(i)}^T \overline{W}_{x_{mj}(i)} \right] = \frac{2}{n_{mj}^2}\sum_{i=1}^{n_{mj}}tr({\overline{\wedge}_{x_{mj}(i)}}^2)
\end{equation}

\begin{equation}\label{var_zsm_app2}
\scriptsize
Var\left[ \frac{1}{n_{mj}} \sum_{i \neq k}^{n_{mj}} \overline{W}_{x_{mj}(i)}^T \overline{W}_{x_{mj}(k)}\right] = \frac{2}{n_{mj}^2}\sum_{i \neq k}^{n_{mj}}tr(\overline{\wedge}_{x_{mj}(i)} \overline{\wedge}_{x(k)})
\end{equation}

The covariance between the second and third term of $Z_{{mj}}$ in (\ref{eq:10}) is 0. Therefore, the variance of $Z_{{mj}}$ is given by adding (\ref{var_zsm_app1}) and (\ref{var_zsm_app2}) which will lead to (\ref{eq:13}).

\section{Proof of Lemma 2}
\label{app_lemma_2}
The cluster compactness of cluster $C_{mj}$ in (\ref{eq:9}) can be written in the following form
\begin{equation}
\begin{split}
y_{{mj}}=\norm{{{x_{mj}} -  x_{mj}}B_{mj}}_2^2  
\end{split}
\end{equation}
where $B_{mj}$ was defined in (\ref{BBB}).

Replacing $x_{mj}(i)$ with ,   $x_{mj}(i) = \overline{c}_{x_{mj}(i)} + \overline{W}_{x_{mj}(i)}$ in (\ref{eq:dd_mod}), we have
\begin{equation}
Y_{{mj}}  = \norm{ \overline{c}_{x_{mj}(i)}A_{mj} + \overline{W}_{x_{mj}(i)}A_{mj} }_2^2  
\end{equation}
\begin{equation} \label{app:3}
\scriptsize
\begin{split}
Y_{{mj}}  =  \norm{ \overline{W}_{x_{mj}(i)}  A_{mj}}_2^2 +
 \frac{2(n_{mj} - 1)}{n_{mj}}\sum_{i=1}^{n_{mj}} \overline{c}_{x_{mj}(i)}^T  \overline{W}_{x_{mj}(i)}\\
 - \frac{2}{n_{mj}} \sum_{i \neq k}^{n_{mj}} \left( \overline{c}_{x_{mj}(i)}^T  \overline{W}_{x_{mj}(k)} + \overline{c}_{x_{mj}(k)} ^T\overline{W}_{x_{mj}(i)} \right) + \norm{\Delta_{{mj}} }_2^2 
\end{split}
\end{equation}
 where 
\begin{equation} \label{app:33}
\begin{split}
\norm{\overline{W}_{x_{mj}(i)}A_{mj}}_2^2  =  \frac{n_{mj}-1}{n_{mj}}\sum_{i=1}^{n_{mj}}\overline{W}_{x_{mj}(i)}^T \overline{W}_{x_{mj}(i)} + \\
 \frac{1}{n_{mj}}\sum_{i \neq k}^{n_{mj}} \overline{W}_{x_{mj}(i)}^T \overline{W}_{x_{mj}(i)}
\end{split}
\end{equation}

Combining (\ref{app:3}) and (\ref{app:3}), $Y_{{mj}}$ can be expressed as follows:
\begin{equation}\label{eq:22}
\begin{split}
Y_{{mj}}  =  \norm{ \Delta_{{mj}} }_2^2 +\frac{n_{mj}-1}{n_{mj}} \sum_{i=1}^{n_{mj}}\overline{W}_{x_{mj}(i)}^T\overline{W}_{x_{mj}(i)}  - \\ \frac{1}{n_{mj}} \sum_{i \neq k}^{n_{mj}}\overline{W}_{x_{mj}(i)}^T\overline{W}_{x_{mj}(k)} +2 \sum_{i=1}^{n_{mj}} \overline{c}_{x_{mj}(i)}^T \overline{W}_{x_{mj}(k)} - \\ \frac{2}{n_{mj}} \sum_{f=1}^{d}\left( \sum_{i=1}^{n_{mj}} \overline{W}_{x_{mj}(k)}(f)  \sum_{i=1}^{n_{mj}}\overline{c}_{x_{mj}(i)}(f) \right)
\end{split}
\end{equation}
where $\overline{W}_{x_{mj}(i)}(f)$ and  $\overline{c}_{x_{mj}(i)}(f)$ are $f$ component of these vectors of length $d$. 

The expected value of $\overline{W}_{x_{mj}(i)}$ is zero, making expectations of terms 4 and 5 of (\ref{eq:22}) to be also zero. The 3rd term's expectation is also zero due to the independence $W$s. Similar to (\ref{exp_inside}), the expectation of the second term can be calculated as follows
\begin{equation}\label{app_ysm_exp}
\begin{split}
E\left[ \frac{n_{mj}-1}{n_{mj}}\sum_{i=1}^{n_{mj}}\overline{W}_{x_{mj}(i)}^T \overline{W}_{x_{mj}(i)} \right] = \frac{n_{mj}-1}{n_{mj}}\sum_{i=1}^{n_{mj}}tr(\overline{\wedge}_{x_{mj}(i)})
\end{split}
\end{equation}
Therefore, combining expectations of all the terms on (\ref{eq:22}), the expected value of $Y_{{mj}}$ is given by (\ref{eq:23}).

For calculation of variance, the first term of equation (\ref{eq:22}) is a constant thus, it has a zero variance. Variances of the other four terms are:
\begin{equation} \label{app:4}
\scriptsize
Var[\frac{n_{mj}-1}{n_{mj}}\sum_{i=1}^{n_{mj}}\overline{W}_{x_{mj}(i)}^T \overline{W}_{x_{mj}(i)}] = 
 \frac{2(n_{mj}-1)^2}{n_{mj}^2}\sum_{i=1}^{n_{mj}}tr((\overline{\wedge}_{x_{mj}(i)})^2)
\end{equation}
\begin{equation}
\scriptsize
Var[ \frac{1}{n_{mj}}\sum_{i \neq k}^{n_{mj}} \overline{W}_{x_{mj}(i)}^T \overline{W}_{x_{mj}(i)}] = \frac{2}{n_{mj}^2}\sum_{i \neq k}^{n_{mj}}tr(\overline{\wedge}_{x_{mj}(i)} \overline{\wedge}_{x_{mj}(i)})
\end{equation}

\begin{equation}
\scriptsize
Var[2\sum_{i=1}^{n_{mj}} \overline{c}_{x_{mj}(i)}^T\overline{W}_{x_{mj}(i)}] = \frac{4}{n_{mj} \ d} \sum_{i=1}^{n_{mj}}  tr(\overline{\wedge}_{x_{mj}(i)}) \sum_{i=1}^{n_{mj}} \overline{c}_{x_{mj}(i)}^T\overline{c}_{x_{mj}(i)}
\end{equation}

\begin{equation}
\scriptsize
\begin{split}
Var[ - \frac{2}{n_{mj}} \sum_{f=1}^{d}\left( \sum_{i=1}^{n_{mj}} \overline{W}_{x_{mj}(i)}(f)  \sum_{i=1}^{n_{mj}}\overline{c}_{x_{mj}(i)}(f) \right)] = \\ \frac{4}{n_{mj}^2 \ d} \sum_{i=1}^{n_{mj}} tr(\overline{\wedge}_{x_{mj}(i)})  \left(\sum_{i=1}^{n_{mj}} \overline{c}_{x_{mj}(i)}^T\overline{c}_{x_{mj}(i)} + \sum_{i \neq k}^{n_{mj}}\overline{c}_{x_{mj}(i)}^T\overline{c}_{x_{mj}(k)} \right)
\end{split}
\end{equation}

The covariance between terms 4 and 5 of equation (\ref{eq:22}) is:
\begin{equation} \label{app:5}
\scriptsize
\begin{split}
2 Cov \left[ 2 \sum_{i=1}^{n_{mj}} \overline{c}_{x_{mj}(i)}^T\overline{W}_{x_{mj}(i)}  - \frac{2}{n_{mj}} \sum_{f=1}^{d}\left( \sum_{i=1}^{n_{mj}} \overline{W}_{x_{mj}(i)}(f)  \sum_{i=1}^{n_{mj}}\overline{c}_{x_{mj}(i)}(f) \right) \right] \\
=\frac{-8}{n_{mj}^2 \ d }\sum_{i=1}^{n_{mj}}tr(\overline{\wedge}_{x_{mj}(i)}) \left(\sum_{i=1}^{n_{mj}} \overline{c}_{x_{mj}(i)}^T\overline{c}_{x_{mj}(i)} + \sum_{i \neq k}^{n_{mj}}\overline{c}_{x_{mj}(i)}^T\overline{c}_{x_{mj}(k)} \right)
\end{split}
\end{equation}

while the rest of the covariances between terms 2-5 of (\ref{eq:22}) are calculated to be zero. Also we have
\begin{equation} \label{app:6}
\scriptsize
\begin{split}
\norm{\Delta_{{mj}} }_2^2 =\norm{\overline{c}_{x_{mj}} A_{mj}}_2^2  =\;\;\;\;\;\;\;\;\;\;\;\;\; \\
\frac{n_{mj} - 1}{n_{mj}}\sum_{i=1}^{n_{mj}}\overline{c}_{x_{mj}(i)}^T\overline{c}_{x_{mj}(i)} - \frac{1}{n_{mj}} \sum_{i \neq k}^{n_{mj}}\overline{c}_{x_{mj}(i)}^T\overline{c}_{x_{mj}(k)} 
\end{split}
\end{equation}

Adding (\ref{app:4})-(\ref{app:5}) and using (\ref{app:6}), the variance of $Y_{{m_j}}$ is given by (\ref{eq:24}).

\section{Upperbound of $\norm{\Delta_{{mj}}}_2^2$ from $y_{{mj}}$ }
\label{app:delta_upp}

To find the upper bound of $  \norm{\overline{c}_{x_{mj}}  A_{mj}}_2^2 $, we should solve for the following inequality:
\begin{equation}\label{eq:31}
E[Y_{{mj}}] - \alpha_{{mj}} \sqrt{var[Y_{{mj}}]} \leq y_{{mj}}
\end{equation}

Using (\ref{eq:23}) and (\ref{eq:24}) in (\ref{eq:31}), and denoting $g_{mj}$ as $g_{mj} = \frac{(n_{mj} - 1)}{n_{mj}} \sum_{i=1}^{n_{mj}}tr(\overline{\wedge}_{x_{mj}(i)}) $, we have in footnote
\footnote{
\begin{equation}
  \norm{\Delta_{{mj}}}_2^2 + g_{mj} - y_{{mj}}   \leq   
  \alpha_{n_{mj}} \sqrt{\frac{2(n_{mj}-1)^2}{n_{mj}^2}\sum_{i=1}^{n_{mj}}tr((\overline{\wedge}_{x_{mj}(i)})^2) + \frac{2}{n_{mj}^2}\sum_{i \neq k}^{n_{mj}}tr(\overline{\wedge}_{x_{mj}(i)} \overline{\wedge}_{x_{mj}(k)}) +   \norm{\Delta_{{mj}}}_2^2\frac{4}{d \times n_{mj}}\sum_{i=1}^{n_{mj}}tr(\overline{\wedge}_{x_{mj}(i)})}\nonumber\;\;\;\;\;(60)
\end{equation}}.
To solve for the boundary, (i.e. when two sides of the inequality are equal), we square both sides of the equation and consequently solve the equation in footnote
\footnote{$(\norm{\Delta_{{mj}}}_2^2)^2 +\norm{\Delta_{{mj}}}_2^2\left(2(g_{mj} - y_{{mj}}) -  \alpha_{n_{mj}}^2 \frac{4}{d \times n_{mj}}\sum_{i=1}^{n_{mj}}tr(\overline{\wedge}_{x_{mj}(i)} )\right) +  $
\begin{equation}
 \;\;\;\;\;\;\;\;\;\;\; \left[ (g_{mj} - y_{{mj}})^2  -\alpha_{n_{mj}}^2 \left(\frac{2(n_{mj}-1)^2}{n_{mj}^2}\sum_{i=1}^{n_{mj}}tr((\overline{\wedge}_{x_{mj}(i)})^2) + \frac{2}{n_{mj}^2}\sum_{i \neq k}^{n_{mj}}tr(\overline{\wedge}_{x_{mj}(i)}  \overline{\wedge}_{x_{mj}(k)} ) \right) \right]=0 \nonumber \;\;\;\;\;\;\;\;\;\;\;\; (61)
\end{equation}}
that is a quadratic equation in terms of $\norm{\Delta_{{mj}}}_2^2$. Solving for the roots of $\norm{\Delta_{{mj}}}_2^2$ can give us both the lower-bound and the upper-bound, denoted by $\overline{\norm{\Delta_{{mj}}}_2^2}$ and $\underline{\norm{\Delta_{{mj}}}_2^2}$, respectively. The higher root of $\norm{\Delta_{{mj}}}_2^2$ is considered its upper bound. Note that the roots $\norm{\Delta_{{mj}}}_2^2$ are only functions of the observed data compactness $y_{m}$, $\alpha_{{mj}}$, and the eigenvalues of the covariance of scatter factor $\overline{W}_{x_{mj}}$, $\overline{\wedge}_{x_{mj}}$.


\begin{thebibliography}{99}
\bibitem{mintro1}
A.K. Jain,
``Data Clustering: 50 Years beyond k-means'',
{\em Patter Recognition}, Lett. 31, vol. 8, pp. 651-666, 2010.


\bibitem{mintro2}
C Aggarwal, C. Reddy,
``Data Clustering: Algorithm and Applications'',
{\em Data Mining and Knowledge Discovery Series},
vol. 31, pp. 684, 2013.

\bibitem{kmeans}
J. MacQueen,
 "Some methods for classification and analysis of multivariate observations." 
 {\em Proceedings of the 5th Berkeley symposium on mathematical statistics and probability.} vol. 1. No. 14. 1967.

\bibitem{ieee1}
X. Dong, P. Frossard, P. Vandergheynst, and N. Nefedov, ``Clustering on Multi-Layer Graphs via Subspace Analysis on Grassmann Manifolds,''
{\em IEEE Transaction on Signal Processing}, vol. 62, no. 4, pp. 905-918, 2014.
\bibitem{ieee2}
B. Yang, X. Fu, and N. D. Sidiropoulos, ``Learning from hidden traits:
Joint factor analysis and latent clustering,''
{\em IEEE Transaction on Signal Processing}, vol. 8, no. 1, pp. 256-269, 2017,
pp. 680-682, Sept. 1985.
\bibitem{ieee3}
Y. Wang, Y. Y. Tang, and L. Li, ``Minimum error entropy based sparse
representation for robust subspace clustering,''
{\em IEEE Transaction on Signal Processing}, vol. 63, no. 15, pp. 4010-4021, 2015.
\bibitem{ref:ditc}
P. Sheng and C. Li, ``Distributed information theoretic clustering,''
{\em IEEE Transaction on Signal Processing}, vol. 62, no. 13, pp. 3442-3453, July
2014.

\bibitem{kmcnc}
M. Chiang, B. Mirkin, 
``Intelligent choice of the number of clusters in kmeans clustering: an experimental study with different cluster spread''
{\em J. Classification},
vol. 27, pp. 3-40, 2010.

\bibitem{gap}
 R. Tibshirani, G. Walther, T. Hastie,
``Estimating the number of clusters in a data set via the gap statistic'',
{\em Journal of the Royal Statistical Society},
 B, 63: pp. 411-423, 2001.
\bibitem{sil}
L. Kaufman and P. Rousseeuw.,  {\em Finding groups in data: an introduction
to cluster analysis.}, New York, NY: Wiley, 1990.
\bibitem{CB}
T. Calinski and J. Harabasz, ``A dendrite method for cluster analysis,'' {\em IEEE
Transactions on Pattern Analysis and Machine Intelligence}, vol. PAMI-1, no. 2, pp. 224-227, April 1979.
\bibitem{DB}
D. L. Davies and D. W. Bouldin. ``A cluster separation measure,''
{\em IEEE Transactions on Pattern Analysis and Machine Intelligence}, vol. PAMI-1, no. 2, pp. 224-227, 1979.

\newpage 

\bibitem{vi1}
N. Yousri, M. Kamek M. Ismail.
``A novel validity measure for clusters of arbitrary shapes and densities''
{\em 19th International Conference on Patter Recognition}
pp. 1-4, 2008.

\bibitem{vi2}
K. Kryszczuk, P. Hurley.
``Estimation of the number of clusters using multiple clustering validity indices'' in 
{\em Multiple Classifier Systems} N. Gayar, J. Kittler, F. Roli, Eds.
Berlin, Heidelberg: Springer, vol. 5997, pp. 114-123, 2010.

\bibitem{xmeans}
D. Pelleg and A. W. Moore, ``X-means: Extending k-means with
efficient estimation of the number of clusters,''
{\em Proceedings of the
Seventeenth International conference on Machine Learning}, pp. 727 -
734, May 2000.
\bibitem{gmeans}
G. Hamerly and Ch. Elkan,
``Learning the K in K-Means,''
 {\em Advances in neural information processing systems},  pp. 281-288,  2004.


\bibitem{dip-means}
A. L. A. Kalegeratos,
``Dip-means: an incremental clustering method
for estimating the number of clusters.''
{\em NIPS’12 Proceeding of the 25th International Conference on Neural Information Processing Systems},
pp. 2393-2401, 2012.

\bibitem{dbscan}
M. Ester, H. Kriegel, J. Sander, and X. Xu,
``A density-based algorithm
for discovering clustering in large spatial databases with noise,''
{\em Proceedings of the Second International Conference on Knowledge Discovery and Data Mining}, pp. 226-231, 1996.
\bibitem{mean_shift}
Y. Cheng,
``Mean shift, mode seeking, and clustering.''
{\em IEEE Transactions on Pattern Analysis and Machine Intelligence}, vol. 17, no. 8, pp. 790-799, August 1995.
\bibitem{optics}
M. Ankerst, M. M. Breuning, H. Kriegel, and J. Sander,
``Optics: Ordering points to identity the clustering structure.'' {\em Proceedings of the 1999 ACM SIGMOD international conference on Management of data}, pp. 49-60, 1999.

\bibitem{fcm}
J. Bezdek, R. Ehrlich, W. Full,
``FCM: The Fuzzy c-Means Clustering algorithm'',
{\em Computers and Geosciences} vol. 10, No. 2-3, pp. 191-203, 1984.

\bibitem{BFC}
T. C. Glenn, A. Zare, and P. D. Gader,
``Bayesian fuzzy clustering,''
\emph{IEEE Transactions on Fuzzy Systems,},
vol. 23, no. 5, pp. 1545-1561, 2015.

\bibitem{MC}
J. Liang, L. Bai, C. Dang, and F. Cao, `The k-means-type algorithms
versus imbalanced data distributions,'' {\em IEEE Transactions on Fuzzy
Systems,} vol. 20, no. 4, pp. 728-745, August 2012.
\bibitem{NML}
S. Hirai and K. Yamanishi,
``Efficient computation of normalized maximum likelihood code for gaussian mixture models with its applications
to clustering,'' {\em IEEE Transactions on Information Theory} vol. 59, no. 11, pp. 7718-7727, 2013.
\bibitem{ref:9}
M. Shahbaba and S. Beheshti, ``MACE-means clustering,'' \emph{
Signal Processing}, vol. 105, pp. 216-225, 2014.
\bibitem{SVD}
 S. Beheshti and S. Sedghizadeh, ``Number of Source Signal Estimation by Mean Squared Eigenvalue Error (MSEE),''
\emph{IEEE Transaction on Signal Processing}, vol. 66, no. 21, pp. 5694-5704, 2018.

\bibitem{ref:8}
S. Beheshti and M. Dahleh, ``Noisy data and impulse response estimation,''
\emph{IEEE Transaction on Signal Processing}, vol. 58, no. 2, pp. 510-521, February 2010.

\bibitem{4}
M. Filipponea, F. Camastrab, F. Masullia, S. Rovettaa.
``A survey of kernel and spectral methods for clustering",
\emph{Pattern Recognition},
vol. 41, pp. 176-190, 2008.

\bibitem{mod_kmeans_p}
K. A. Abdul Nazeer, M. P. Sebastian.
``Improving the Accuracy and Efficiency of the k-means Clustering Algorithm",
\emph{Proceedings of the World Congress on Engineering}, 
vol 1, 2009

\bibitem{kkmace_p}
F. Rahman, S. Beheshti.
``Kernel k-MACE Clustering",
\textit{52nd Asilomar Conference on Signals, Systems, and Computers}, 
2018



\bibitem{ref:8.1}
S. Beheshti and M. Dahleh,
``A new information-theoretic approach to signal denoising and best basis selection,''
{\em IEEE Transaction on Signal Processing}, vol. 53, no. 10, pp. 3613-3624, October 2005.

\bibitem{ari}
K. Yeung, W. Ruzzo, 
"Details of the adjusted rand index and clustering algorithms: Supplement to the paper 'an empirical study on principal component analysis for clustering gene expression data’",
{\em Bioinform}, no. 17, pp. 763-774, 2001.

\bibitem{nvi}
M. Meila, "Comparing Clusterings by the Variation of Information",
{\em Proc. Conf. Learning Theory}, 2003.

\bibitem{d31}
A. Gionis, H. Mannila, and P. Tsaparas,
``Clustering aggregation,''
{\em ACM
Transactions on Knowledge Discovery from Data}, vol. 1, no. 1, pp. 1-30, 2007.
\bibitem{s15}
C. Veenman, M. Reinders, and E. Backer,
``A maximum variance cluster algorithm,''
{\em IEEE Transactions Patter Analysis and Machine Learning}, vol. 24, no. 9, pp. 1273-1280, 2002.





\end{thebibliography}

\end{document}